\documentclass[authoryear,1p,nopreprintline,12pt]{elsarticle}
\usepackage{geometry}
\usepackage[english]{babel}
\usepackage[T1]{fontenc}
\usepackage[hidelinks]{hyperref}
\usepackage{amsmath, amssymb}
\usepackage{algorithm}
\usepackage{algpseudocode}
\usepackage{enumitem}
\usepackage{booktabs}
\usepackage{multirow}
\usepackage{mathtools}
\usepackage{tikz}
\usepackage{dsfont}
\usepackage{pgfplots}
\pgfplotsset{compat=1.18}
\usetikzlibrary{arrows.meta,positioning,calc,fit,backgrounds}
\usepackage{graphicx}
\usepackage{subcaption}
\usepackage{float}
\usepgfplotslibrary{groupplots}
\usepackage{xcolor} 
\usepackage{etoolbox}
\usepackage{amsthm}

\usepackage[title]{appendix}

\usepackage{cleveref}
\usepackage{csquotes}
\usepackage{siunitx}
\usepackage{soul} 

\definecolor{mpiblue}{HTML}{33A5C3}
\definecolor{mpigreen}{HTML}{007675}
\definecolor{mpired}{HTML}{78004B}
\definecolor{mpiblack}{HTML}{383C3C}
\definecolor{mpigray}{HTML}{87878D}

\begin{document}

\begin{frontmatter}

\title{Optimization with Region-Reduced ReLU Neural Networks}

\author[MPI_OML_affiliation,Uni_affiliation]{Christoph Plate}
\author[MPI_PSE_affiliation,LUH_EVT_affiliation]{Caroline Ganzer}
\author[Uni_affiliation]{Mirko Hahn}
\author[MPI_PSE_affiliation]{Alexander Klimek}
\author[Uni_affiliation]{Heyuan Liu}
\author[MPI_OML_affiliation,Uni_affiliation]{Sebastian Sager}
\author[MPI_PSE_affiliation,Uni_affiliation_SVT]{Kai Sundmacher}
\author[Uni_affiliation]{Hanna Wilhelm}

\affiliation[MPI_OML_affiliation]{
            organization={Max Planck Institute for Dynamics of Complex Technical Systems, Mathematical Optimization and Machine Learning Group},
            addressline={Sandtorstraße~1}, 
            city={Magdeburg},
            postcode={39106}, 
            state={Saxony-Anhalt},
            country={Germany}}
            
\affiliation[Uni_affiliation]{
            organization={Otto~von~Guericke~University, Faculty of Mathematics, Insitute of Mathematical Optimization},
            addressline={Universitätsplatz~2}, 
            city={Magdeburg},
            postcode={39106}, 
            state={Saxony-Anhalt},
            country={Germany}}            

\affiliation[MPI_PSE_affiliation]{
            organization={Max Planck Institute for Dynamics of Complex Technical Systems, Department of Process Systems Engineering},
            addressline={Sandtorstraße~1}, 
            city={Magdeburg},
            postcode={39106}, 
            state={Saxony-Anhalt},
            country={Germany}}
            
\affiliation[LUH_EVT_affiliation]{
            organization={Leibniz~University~Hannover, Faculty of Mechanical Engineering, Insitute of Technical Combustion, Chair of Energy~Process~Engineering},
            addressline={An~der~Universität~1}, 
            city={Garbsen},
            postcode={30823}, 
            state={Lower~Saxony},
            country={Germany}}            

\affiliation[Uni_affiliation_SVT]{
            organization={Otto~von~Guericke~University, Faculty of Process and Systems Engineering, Insitute of Process Engineering, Chair of Process Systems Engineering},
            addressline={Universitätsplatz~2}, 
            city={Magdeburg},
            postcode={39106}, 
            state={Saxony-Anhalt},
            country={Germany}}

\begin{abstract}
Optimization of mathematical models involving integer decisions and neural networks with ReLU activation (ReLU ANNs) is a challenging task. Nevertheless, such models are an enabling technology in many application domains. A prominent example is superstructure optimization in chemical engineering, where ReLU ANNs are frequently employed as surrogate models for complex nonlinear processes.

We survey recent developments in this area. We argue that in addition to network size and training options of the ANNs, the ReLU activation geometry and the number of linear regions on the domain of interest have a strong impact on computational optimization performance. While standard model compression approaches such as structured pruning reduce network size, they do not explicitly address geometric considerations.  Therefore, we propose a novel \emph{region-reduced} model compression approach that combines the stabilization of unstable neurons and the merging of redundant neurons to reduce the number of linear regions while maintaining predictive accuracy through error compensation. We evaluate our method against standard compression approaches on multiple optimization use cases. First, the two-dimensional peaks function for which we can visualize the activation geometry. Second, on optimization over individual surrogate ReLU ANNs for three chemical processes, and third, on a hybrid superstructure optimization problem that involves the three ReLU ANNs, additional process submodels, and binary variables. The results for the superstructure problem demonstrate the large potential of region reduction with a decrease of 40\% to 50\% in computational time, yielding solutions closer than 1\% to the reference at negligible effort of obtaining the compressed model.
\end{abstract}

\begin{keyword}
ReLU neural network \sep Mixed-integer linear programming \sep 
Model compression \sep Superstructure optimization
\end{keyword}
\end{frontmatter}

\section{Introduction}\label{sec:intro}
We consider feed-forward artificial neural networks (ANNs) of the following form,
\begin{align}\label{eq:ANN}
	x^{(\ell)} = \sigma^{(\ell)}\left( W^{(\ell)} x^{(\ell-1)} + b^{(\ell)}  \right), \quad \ell \in [L].
\end{align}
Here and in the following the superscript $(\ell)$ will refer to one out of $L$ layers of the ANN. 
The input of the ANN is denoted by $x^{(0)} = x \in \mathbb{R}^{n_0}$ and $W^{(\ell)} \in \mathbb{R}^{n_{\ell}\times n_{\ell-1}}, b^{(\ell)} \in \mathbb{R}^{n_{\ell}}$ are the weights and biases of layer $\ell$, respectively.
On all layers an activation function $\sigma^{(\ell)}$ is applied to an affine transformation of the layer input. We consider the special case
\begin{align*}
    \sigma^{(\ell)}(x) = \max\{0,x\} \quad \forall \; \ell \in [L-1], \qquad \sigma^{(L)}(x) = x
\end{align*}
of the ReLU activation function on all but the last layer and call the resulting model \textit{ReLU ANN}.

Embedding such pre-trained ANNs into optimization problems has become ubiquitous in many areas of application. Examples include model-predictive control with learned dynamics~\citep{Salzmann2023}, operations research~\citep{Dixit2025} with data-driven cost/revenue surrogates, and superstructure optimization in process systems engineering~\citep{Klimek2026}. The universal approximation property~\citep{Cybenko1989} makes ANNs a powerful tool, particularly in cases where data is widely available and rigorous models are either unavailable or too costly to consider. 

Different activation functions $\sigma$ have been studied in the context of optimization with embedded ANNs.
Popular choices are smooth activation functions that allow the application of gradient-based, convexification, and spatial Branch and Bound methods~\citep{Schweidtmann2019a}.
The ReLU activation function $\textrm{ReLU}(x) = \max\{0,x\}$ in contrast has a combinatorial nature, due to the maximum operator~\citep{Fischetti2018}. With the big-$M$ method a neuron with ReLU activation function can be formulated using binary decision variables $z_i^{(\ell)}$ that indicate if a neuron is ``active'', i.e., if $\max\{0,x\} = x$. This formulation uses valid pre-activation bounds $L^{(\ell)} \leq  W^{(\ell)} x^{(\ell-1)} + b^{(\ell)} \leq U^{(\ell)}$ and linear inequalities,
\begin{align}\label{eq:bigM}
\begin{split}
    x_i^{(\ell)} &\ge 0, \\
    x_i^{(\ell)} &\ge  \left(W^{(\ell)} x^{(\ell-1)} + b^{(\ell)}\right)_i, \\
    x_i^{(\ell)} &\le \left(W^{(\ell)} x^{(\ell-1)} + b^{(\ell)}\right)_i - L_i^{(\ell)} \left(1-z_i^{(\ell)}\right), \\
    x_i^{(\ell)} &\le U_i^{(\ell)} z_i^{(\ell)}, \\
    z_i^{(\ell)} & \in \{0,1\},
\end{split}
\end{align}
that are incorporated for all neurons in the ANN.
As a result, using ReLU activation functions is a practical way to embed nonlinear expressions in mixed-integer linear (MILP) or mixed-integer quadratic (MIQP) optimization problems. This makes the problem amenable to global optimization without needing to resort to general (mixed-integer) nonlinear optimization.

In practice, however, even state-of-the-art optimization solvers such as Gurobi~\citep{gurobi} can handle only relatively small-scale ANNs with no more than a few hidden layers. 
There are two major aspects that have a strong impact on the computational time of mixed-integer programming solvers in the context of ReLU ANNs.
First, the tightness of relaxations of the binary variables $z_i^{(\ell)}$ and algorithmic approaches. In particular, the magnitude of the big-$M$ coefficients $L^{(\ell)}, U^{(\ell)} \in \mathbb{R}^{n_{\ell}}$ in \eqref{eq:bigM} may grow exponentially in the layer number $\ell$, resulting in weak relaxations and a strong increase in overall computational time. 
Second, the number of binary variables $z_i^{(\ell)}$ is proportional to the total number of hidden neurons in the network. While deep learning is desirable from an ANN accuracy point of view, the runtime of mixed-integer programming is typically exponential in the number of binary variables. 

Several methods have been proposed to reduce the effort of solving optimization problems with embedded ReLU ANNs. We survey previous work before we outline our contribution, a novel geometrically-inspired compression method. For a recent review of related software packages we refer to \citet{Plate2026}.

\subsection{Optimization over piece-wise linear submodels}

An early detailed study of optimization of models including piecewise linear geometry was given in~\citet{Joseph-Duran2014}. There, the functions $\max\{0, x\}$ modeled the overflow of sewage water. It was shown that MILP reformulations or smoothing may be inferior to tailored algorithms. However, the case of embedding ReLU ANNs is more general in the sense that the objective function of the optimization problem does not necessarily point in the direction of $0$, as in the case of water overflow which shall be avoided.
Thus, for ReLU ANNs the big-$M$ formulation is still the most popular approach in the literature \citep{Fischetti2018,Xiao2019,Tjeng2019,Yang2022}. 

MILP formulations of ReLU neurons with tighter relaxations compared to the standard big-$M$ formulation and cuts derived from those formulations have been investigated \citep{Anderson2020,Tsay2021}. Also tight formulations of piecewise-linear functions in MILPs based on a difference-of-convex functions representation were proposed \citep{Ploussard2025}.
Different bound-tightening methods \citep{Grimstad2019,Badilla2023} are popular in mixed-integer programming. Moreover, several exact transformations of ReLU ANNs have been proposed that yield functionally equivalent ANNs with beneficial properties for optimization, e.g., lower big-$M$ coefficients and less binary variables \citep{Kumar2019,Plate2026}. Recently, also input-convex neural networks were investigated as surrogates in optimization problems \citep{Liu2026}, and the authors reported gains in optimization time compared to standard neural networks in cases in which the modeled functions are convex.

\subsection{Obtaining and tightening big-M coefficients}

Tractability of optimization problems with embedded ReLU ANNs hinge on the magnitude of big-$M$ coefficients occurring in formulation \eqref{eq:bigM}. These bounds can be determined via interval arithmetic by propagating known input bounds $L^{(0)}, U^{(0)}$ with $x \in [L^{(0)}, U^{(0)}]$ through the network, i.e.,

\begin{align}\label{eq:interval_arithmetic}
    \begin{split}
        L^{(\ell)}_i &= \sum_{j=1}^{n_{k-1}} \min \left\{ W_{i,j}^{(\ell)} L_j^{(\ell-1)}, W_{i,j}^{(\ell)} U_j^{(\ell-1)} \right\} + b^{(\ell)}_i, \quad \ell \in [L], \ i \in [n_\ell],  \\
        U^{(\ell)}_i &= \sum_{j=1}^{n_{k-1}} \max \left\{ W_{i,j}^{(\ell)} L_j^{(\ell-1)}, W_{i,j}^{(\ell)} U_j^{(\ell-1)} \right\} + b^{(\ell)}_i, \quad \ell \in [L], \ i \in [n_\ell].
    \end{split}
\end{align}

Despite being a simple method, existing dependencies in the activation status of neurons are not taken into account which leads to loose bounds.
Optimization-based bound-tightening (OBBT) can exploit these dependencies to obtain tighter bounds based on one additional minimization and maximization of each neuron's pre-activation value \citep{Grimstad2019,Badilla2023}. To reduce computational effort, typically the relaxed problem is solved. E.g., the LP to obtain $L^{(\ell)}_k$ reads 

\begin{equation}\label{prob:obbt}
    \begin{aligned}
    \min_{x,z}\ & W^{(\ell)}_k x^{(\ell-1)} + b^{(\ell)}_k \\ 
    \textrm{s.t.}\ & \begin{alignedat}[t]{3}
            x^{(j)}_i & \geq 0, \quad && j \in [\ell],\, i \in [n_j], \\
            x^{(j)}_i & \geq W^{(j)}_i x^{(j-1)} + b^{(j)}_i, \quad && j \in [\ell],\, i \in [n_j], \\
            x^{(j)}_i & \leq \left(W^{(j)}_i x^{(j-1)} + b^{(j)}_i\right) - L^{(j)}_i \left(1-z^{(j)}_i\right), \quad && j \in [\ell],\, i \in [n_j], \\
            x^{(j)}_i & \leq U^{(j)}_i z^{(j)}_i, \quad && j \in [\ell],\, i \in [n_j], \\
            L^{(0)}_i & \leq x^{(0)}_i  \leq U^{(0)}_i,        \quad && i \in [n_x],\\
            z^{(j)}_i & \in [0,1], \quad && j \in [\ell],\, i \in [n_j].
        \end{alignedat}
    \end{aligned}
\end{equation}

\subsection{The ReLU ANN training phase}

The ANN's properties are formed during training. By selecting a suitable optimization algorithm, objective function and hyperparameters, these properties can be influenced. Hence, from a practical point of view when considering optimization with embedded ReLU ANNs it makes a difference whether the training process of the considered ReLU ANN can still be influenced, or if one is presented with a fixed network that has already been trained. Previous work has shown that the choice of training hyperparameters has a strong impact on the computational tractability of subsequent optimization problems.
Notably, standard $\ell_1$ weight regularization applied during training was shown to be capable of increasing ReLU stability and reducing big-$M$ coefficients \citep{Xiao2019,Plate2026}. Moreover, tailored regularization terms have been proposed to increase tractability of downstream MILP problems by promoting ReLU stability \citep{Xiao2019} or by penalizing either the magnitudes of big-$M$ coefficients or the LP relaxation gap during training \citep{Tsay2026}. 
In contrast, the use of redundancy-enhancing methods like dropout has been found to increase the difficulty of solving the resulting optimization problems in certain settings \citep{Plate2026}. 
Here, we focus on the situation in which the original ReLU ANN weights are given and focus on model compression techniques.

\subsection{Compressing ReLU ANNs}

If the training phase cannot be influenced but the network is too large to be directly incorporated into an optimization problem, various model compression techniques can be employed.
Pruning is an established method to reduce the size of ANNs. Dating back to \citet{Janowsky1989}, the most commonly used pruning method is based on the magnitude of the weights, and is called \emph{Magnitude Pruning}. The basic idea is to remove the weights with the smallest absolute value, usually followed by a fine-tuning of the remaining parameters to recover the desired level of accuracy. Pruning methods can be generally divided into two groups: \emph{Unstructured} pruning methods remove single connections, i.e., single entries in the weight matrices $W^{(\ell)}$, which leads to sparser weight matrices, whereas \emph{structured} pruning removes parameters in groups, e.g., complete neurons \citep{Blalock2020}.
Several papers studied pruning methods and their effect on approximation qualities \citep{Han2015,Suzuki2020a}. In recent years, with the emerging trend of addressing optimization problems with embedded neural networks, pruning methods have also gained importance as a pre-processing step to facilitate and accelerate their solution. In~\citet{Cacciola2024}, structured pruning was applied to facilitate the solution of verification problems, i.e., determining whether an adversarial perturbation exist that incorrectly changes the prediction of the classifier given by the ANN. A heuristic to solve these problems was introduced in \citet{Pham2026}, where it was observed that even coarse surrogate models obtained via unstructured pruning without a fine-tuning step can speed up the solution of verification problems. Although there exist more elaborate methods, e.g., iterative pruning methods \citep{Lastrucci2026}, standard Magnitude Pruning, i.e., pruning weights with small absolute values or neurons with smallest row-wise norms in the weight matrices, remains a popular method due to its efficiency and efficacy \citep{Pham2026}. 

Similar to pruning, knowledge distillation is a compression technique for neural networks where a smaller student model is trained from scratch to replicate the behavior of a larger, pre-trained teacher model \citep{Ba2014,Hinton2015}. In the offline distillation setting, the teacher model is pre-trained on training data. Afterwards, this trained model is used to extract training labels which provide targets for the student to be trained on, enabling it to learn the teacher's behavior. 
\citet{Gou2021} survey knowledge distillation, especially for classification tasks. Regression is considered in \citep{Kang2021,Zhou2023}.

\subsection{Stable active or inactive neurons}

ReLU-stability refers to whether the activation pattern of ReLU neurons remains unchanged under perturbation in the inputs $x$ \citep{Xiao2019,Kumar2019}. Extreme cases -- neurons that are constantly active or inactive for all values of $x$ -- do not need to be modeled with binary variables. They can be equipped with the identity instead of the ReLU activation or be removed from the network, respectively.
For neurons that are either active or inactive for almost all input vectors $x$, modeling them without binary variables implies an approximation and thus a potential loss in accuracy. To understand and leverage the trade-off between reduced accuracy and computational speed-up, a more detailed analysis is necessary. For this, the concept of \textit{linear regions} is helpful which we formalize next.

A feed-forward ReLU neural network with $L$ layers of $\{n_\ell\}_{\ell=1}^L$ neurons per layer with a total number of $N=\sum_{\ell=1}^{L-1} n_\ell$ hidden neurons implements a continuous piecewise-affine map ${f:\mathbb{R}^{n_0} \mapsto \mathbb{R}^{n_L}}$ \citep{Grigsby2022}. Let $g^{(\ell)}(x) := W^{(\ell)} x^{(\ell-1)} + b^{(\ell)} $ be the pre-activations in layer $\ell$. We define the vectors $z^{(\ell)}\in\{0,1\}^{n_{\ell}}$ with $z^{(\ell)}_i := \mathds{1}_{\{g^{(\ell)}_i(x)\ge0\}}$ as the \emph{activation pattern} of layer $l \in [L-1]$. Within a given domain of interest $\mathcal D\subset\mathbb{R}^{n_0}$, each feasible activation pattern induces a convex polytope (called a \emph{linear region} in the following)
\begin{align}
\begin{split}
\label{eq:region-polytope}
\mathcal{R}(z;\mathcal D)
=\bigl\{x\in \mathcal D:\; & g_i^{(\ell)}(x)\ge 0 \text{ if } z^{(\ell)}_i=1,\;\; \\ 
 & g_i^{(\ell)}(x) \le 0 \text{ if } z^{(\ell)}_i=0,\ \forall \ell \in [L-1] , i \in [n_\ell])\bigr\},
 \end{split}
\end{align}
on which $f$ is affine. We call the number of distinct nonempty linear regions on the domain $\mathcal D$ the (domain-wise) \emph{number of linear regions}, denoted $\#\mathrm{LR}(f;\mathcal D)$, or $\#\mathrm{LR}$ for short. This definition of linear regions is not unique, alternative definitions have been proposed, e.g., in \citet{Stargalla2026}.

Upper and lower bounds on the number $\#\mathrm{LR}(f;\mathcal D)$  of linear regions have been studied extensively. The first trivial upper bound was derived in \citet{Montufar2014} by associating hidden neurons with ReLU activation patterns as
\begin{align}
    \#\mathrm{LR}(f;\mathcal D) \leq 2^{N}.
\end{align}
This bound was improved in \citet{Montufar2017} to
\begin{align}
    \#\mathrm{LR}(f;\mathcal D) \leq \prod_{\ell=1}^{L-1} \sum_{j=0}^{d_{l}} \binom{n_{\ell}}{j}
\end{align}
with $d_\ell=\min\{n_0,\ldots,n_{\ell-1}\}$. Later, a stronger bound was derived in \citet{Serra2018} with
$J = \{ (j_1,\ldots,j_{L-1}) \in \mathbb{Z}^{L-1}: 0 \leq j_\ell \leq \min\{n_0,n_1-j_1,\ldots,n_{\ell-1}-j_{l-1},n_{\ell}\}, \ \forall \ell \in [L-1] \}$
as
\begin{align}
    \#\mathrm{LR}(f;\mathcal D) &\leq  \sum_{(j_1,\ldots,j_{L-1}) \in J} \prod_{\ell=1}^{L-1} \binom{n_{\ell}}{j_\ell}
\end{align}
by identifying small layers as bottlenecks that impede the full capacity of the network to form linear regions. Still, in practice one often finds that these bounds are typically not reached by trained neural networks \citep{Hanin2019,Serra2018}. In \citet{Gamba2022}, the authors argue that counting linear regions alone may also overestimate a network's effective nonlinearity. They analyze trained networks and show that many linear regions are locally redundant, i.e., they are associated with similar affine functions. Moreover, it has been noted that training options affect the number and size of the linear regions \citep{Zhang2020a}.

The notion of linear regions has been used in different contexts of optimization with embedded ReLU ANNs. First, in \citet{Tong2024a} an optimization heuristic for the special case of optimizing a linear function of a neural network's output over polytopes is proposed that makes use of the network's linear regions. By iteratively solving the LP relaxation within the current linear region and then computing a step into the neighboring linear region until no further improvement is achieved, the authors are able to report better scalability of their heuristic compared to standard MILP solvers for ANNs with increasing size, while still obtaining high-quality solutions. Second, in \citet{Plate2026}, the influence of training options on the number of linear regions and computational complexity of subsequent optimization problems has been investigated.

In summary, the number of linear regions of an ANN is highly relevant for optimization with embedded ReLU ANN. This number corresponds to the feasible binary solutions to \eqref{eq:bigM} for all neurons and thus to the maximum depth of a Branch and Bound tree only for the binary variables corresponding to the ReLU ANN. 
We are not aware of any previous work on explicitly reducing the number of linear regions in an optimization context.

\subsection{Main contributions and outline}

We propose a new and efficient geometric post-training pipeline that reduces the number of linear regions of trained ReLU neural networks while maintaining a high predictive accuracy. It employs a bias-shift routine that moves boundaries of linear regions by adjusting the bias values of specific neurons and a subsequent merging of neurons that have highly correlated activation patterns. Geometrically, the approach reduces the number of linear regions as specified in \eqref{eq:region-polytope}. We shall refer to this as \textit{region-reduction} in the following. 

The paper is structured as follows. In \Cref{sec:method} we develop the new method with its two main ingredients bias-shift and pattern-merge. 
We describe different benchmark applications in \Cref{sec:applications}.
In \Cref{sec:exper} we evaluate our method on these applications, benchmark it against state-of-the-art compression methods, i.e., magnitude pruning and knowledge distillation, and discuss the results. We conclude with a summary.

\newcommand{\relu}{\mathrm{ReLU}}
\DeclarePairedDelimiter{\norm}{\lVert}{\rVert}
\DeclarePairedDelimiter{\inner}{\langle}{\rangle}
\newcommand{\diam}{\mathrm{diam}}
\newcommand{\R}{\mathbb{R}}

\section{Methods and Models} \label{sec:method}

In this section we define the four types of ReLU ANNs that are embedded in the optimization problems and compared for performance in the next section. 
The \textbf{original} model is a ReLU ANN of type \eqref{eq:ANN} with given weights and biases, trained on synthetic data generated by simulating nonlinear functions. However, the same procedure would apply to ReLU ANNs inferred from real-world data. Previous work \citep{Plate2026} shows that optimization performance strongly depends on the regularization during the training process. Therefore, throughout our experiments we consider one version without regularization and one with an $\ell_1$ penalty of the network weights, scaled with a hyperparameter $\lambda$.

\begin{figure}[h]
    \centering
    \includegraphics[width=\linewidth]{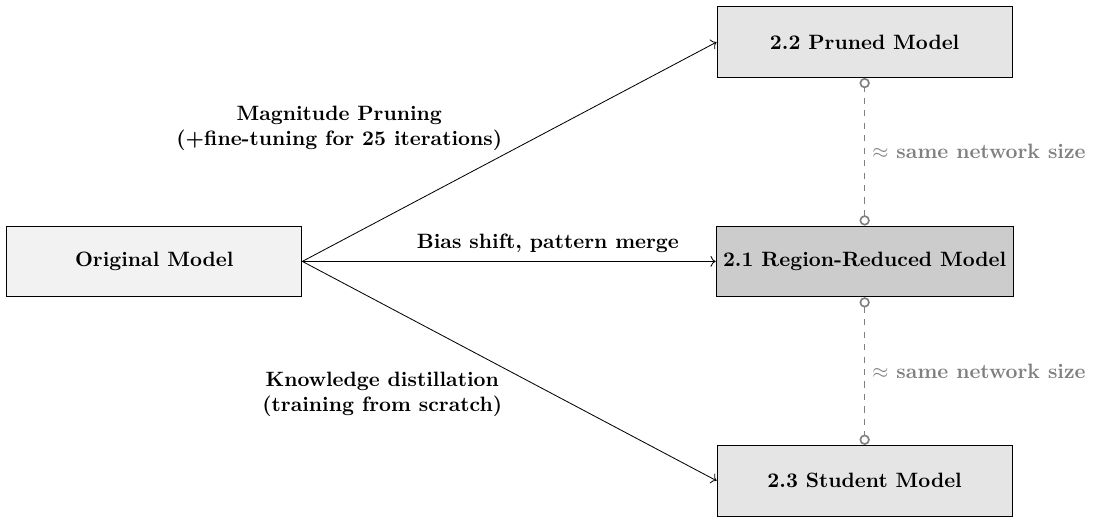}
    \caption{The four model types compared in this study: the original model and three derived compressed models. While the compression ratio of the region-reduced model is varying for each network, we use a fixed compression ratio of 20\% for the pruned and the student model.}
    \label{fig:workflow}
\end{figure}

The other three model classes are derived from the original model as indicated in Figure~\ref{fig:workflow}. The \textbf{region-reduced} model is novel, and obtained via our proposed method outlined in \Cref{subsec:ReRe}. The \textbf{pruned} and \textbf{student} models have already been established in the literature and can thus be used as a benchmark. They are shortly introduced in \Cref{subsec:pruned,subsec:student}. All three approaches are applicable after the training of the original model, i.e., independent of the original training data.

\subsection{Region-Reduced Model} \label{subsec:ReRe}

Our approach for post-processing of ReLU ANNs after their training consists of three steps aimed at reducing the total number of linear regions.  In a first step (\textit{bias shift}), neurons which are almost always off are stabilized and subsequently removed from the network. Then, similar neurons are merged (\textit{pattern merge}), lowering the number of neurons in the network. Lastly, \textit{bias shift} is used again to stabilize neurons that are almost always active. Both methods \textit{bias shift} and \textit{pattern merge} depend on an evaluation of ReLU stability and  rigorous pre-activation bounds $\{L^{(\ell)},U^{(\ell)}\}_{\ell=1}^L$ obtained either via \eqref{eq:interval_arithmetic} or \eqref{prob:obbt} to guide the selection of neurons. Using \eqref{eq:ANN}, we define the
\textit{pre-activation} functions 
$$g^{(\ell)}(x) := W^{(\ell)} x^{(\ell-1)} + b^{(\ell)}$$ 
and the \textit{post-activation} functions 
$$h^{(\ell)}(x) := \relu\left(g^{(\ell)}(x)\right).$$
We identify candidate neurons by calculating statistics of the pre-activation values of all ReLU neurons on a probe batch $X \subset \mathbb{R}^{n_0}$ which is sampled from the input domain $\mathcal D$. Each of the $n_X$ input samples $x^{(k)} \in \mathbb{R}^{n_0}$ is propagated through the ANN.  We define the vector of pre- and post-activation values for the samples of the probe batch ${X}$ for neuron $i$ on layer $\ell$ as 
$$g_i^{(\ell)}(X) := \left( \begin{array}{c} g_i^{(\ell)}(x^{(1)}) \\ \dots \\ g_i^{(\ell)}(x^{(n_X)}) \end{array} \right)$$ 
and
$$h_i^{(\ell)}(X) := \left( \begin{array}{c} h_i^{(\ell)}(x^{(1)}) \\ \dots \\ h_i^{(\ell)}(x^{(n_X)}) \end{array} \right).$$ 
Similarly, with $g^{(\ell)}(X), h^{(\ell)}(X) \in \mathbb{R}^{n_\ell \times n_X}$ we denote the matrices with pre- and post-activation values of neurons in layer $\ell$ for all samples in the probe batch $X$.
To assess activity of individual neurons, we define the sampled activation statistic as
\begin{align}\label{eq:activation_approximated}
a_i^{(\ell)} := \frac{1}{n_X} \sum_{k=1}^{n_X} \mathds{1}_{\{h_i^{(\ell)}(x^{(k)}) > 0\}} 
\end{align}
for $\ell \in {[L-1]}, i \in [n_\ell]$ and call a neuron \emph{almost stable} if $a_i^{(\ell)} \le \varepsilon$ or $a_i^{(\ell)} \ge 1-\varepsilon$. E.g., with $\varepsilon = 0.01$ a neuron is almost stable if its activation $z_i^{(\ell)}$ is the same for at least $99\%$ of the samples in the probe batch $X$. 

\subsubsection{Bias Shift} \label{subsec:bs}
In a first step, we adjust the bias vectors $b^{(\ell)}, \ \ell \in {[L-1]}$ of almost stable neurons. Geometrically, this corresponds to a translation of the (bent) hyperplanes induced by the specific neurons.

To control the magnitude of perturbations, we impose an upper limit on the allowed bias shifts. This limit is chosen adaptively per neuron, i.e., as a fraction $\rho \in [0,1]$ of the range of observed pre-activation values on the probe batch $X$. Mathematically, this upper limit is defined as 

\begin{equation}
    \delta_{\max,i}^{(\ell)} = \rho \cdot \left( \max_{x \in X} g^{(\ell)}_i(x) -  \min_{x \in X} g^{(\ell)}_i(x)  \right).
\end{equation}

The method can be used in two different settings: either to permanently deactivate neurons that are almost always inactive on the probe batch or to permanently activate neurons that are almost always active. The necessary bias shift to deactivate or activate such neurons is computed as 

\begin{equation}\label{eq:delta}
    \Delta^{(\ell)}_i = 
\begin{cases}
    s^{(\ell)}_i \left(m + \min_{x \in X} s^{(\ell)}_i g_i^{(\ell)}(x)\right), & \text{if } |\min_{x \in X} s^{(\ell)}_i  g_i^{(\ell)}(x) | \leq  \delta_{\max,i}^{(\ell)} \\
    0, & \text{else}
\end{cases}
\end{equation}
with $m$ being a margin added for safety and $s^{(\ell)}_i \in \{-1,1\}$ being the rounded average activation status of the neuron indicating the direction of the shift, i.e.,
\begin{equation}
s^{(\ell)}_i = 2 \left\lfloor a_i^{(\ell)}  \right\rceil - 1.
\end{equation}

After permorming these operations for all neurons $i \in [n_{\ell}]$, the bias vector on the current layer $\ell$ is updated as $b^{(\ell)}+\Delta^{(\ell)}$. Before continuing the procedure on the next layer with ReLU activation, the propagated error can be reduced by modifying the biases of this next layer. We estimate the average deviation in the post-activations of the current layer introduced by the bias shifts as
\begin{align}
    \Delta h^{(\ell)} = \frac{1}{n_X} \sum_{k=1}^{n_X} \left( \relu\left(g^{(\ell)}(x^{(k)}) + \Delta^{(\ell)}\right) - \relu\left(g^{(\ell)}(x^{(k)})\right)  \right) .
\end{align}
Multiplying this estimated error $\Delta h^{(\ell)}$ with the outgoing weights and subtracting the result from the subsequent layers' biases yields an updated bias vector that compensates the introduced error, i.e.,
\begin{align}
    \begin{split}
        b^{(\ell+1)} \gets b^{(\ell+1)} -  W^{(\ell+1)} \Delta h^{(\ell)}.
    \end{split}
\end{align}

After applying error correction to compensate for the bias adjustments, neurons that were stabilized to remain inactive are pruned from the neural network. This step reduces the model's parameter count while preserving its functional behavior. The information on neurons that are stabilized to remain active is saved and returned at the end of the overall algorithm, such that the corresponding binary variables of a subsequent MILP embedding can be fixed. Formally, one could also change the activation function of the stabilized neurons to the identity. In practice however, the workaround with fixed binary variables is necessary because the frameworks for embedding neural networks in optimization problems, e.g., OMLT \citep{Ceccon2022}, expect layers with one activation function for all neurons. \Cref{alg:biasshift} summarizes the outlined procedure.

\begin{algorithm}[H]
\caption{\textsc{BiasShift}}
\label{alg:biasshift}
\begin{algorithmic}[1]
\Require Weights {$\{W^{(\ell)}, b^{(\ell)}\}_{\ell=1}^{L}$}, bounds $\{L^{(\ell)}, U^{(\ell)}\}_{\ell=1}^{L}$, data {$X=\{x^{(k)}\}_{k=1}^{n_X}$}, direction {$d \in \{\text{off},\text{on} \}$}
\Require Margin $m = 10^{-6}$, bias shift fraction  $\rho = 0.3$, stability threshold $\varepsilon = 0.1$.
\State Compute $\{g^{(\ell)}(X)\}_{\ell=1}^L$ and  $\{h^{(\ell)}(X)\}_{\ell=1}^L$ \Comment{Initial forward pass}
\State $\mathcal{S} \gets \emptyset$ \Comment{Set of stable and stabilized neurons} 
\For{$\ell =1$ \textbf{to} $L-1$}
  \For{$j = 1$ \textbf{to} $n_\ell$}
    \State $\Delta^{(\ell)}_j \gets 0$ \Comment{Initialize with zero}
    \If{($d = \text{off}$ and $U_j^{(\ell)} < 0$) or ($d = \text{on}$ and $L_j^{(\ell)} > 0$)} 
        \State $\mathcal{S} \gets \mathcal{S} \cup \{(\ell,j)\}$
        \State continue
    \EndIf
    \State $a_j^{(\ell)} = \frac{1}{n_X} \sum_{k=1}^{n_X} \mathds{1}_{\{h_j^{(\ell)}(x^{(k)}) \geq 0\}}$ \Comment{Sampled activation statistic}
    \If{($d = \text{off}$ and $a_j^{(\ell)} \le \varepsilon$) or ($d = \text{on}$ and $a_j^{(\ell)} \ge 1-\varepsilon$)} 
        \State $\delta_{\max,j}^{(\ell)} \gets \rho \cdot \left( \max_{x \in X} g^{(\ell)}_j(x) -  \min_{x \in X} g^{(\ell)}_j(x)  \right)$ \Comment{Upper limit}
        \State $s_j^{(\ell)} \gets 2 \left\lfloor a_j^{(\ell)}  \right\rceil - 1$ 
        \If{$|\min_{x \in X} s^{(\ell)}_j \cdot g_j^{(\ell)}(x)| \leq \delta_{\max,j}^{(\ell)}$} 
            \State $\Delta^{(\ell)}_j \gets  s^{(\ell)}_j \cdot\left(m- \min_{x \in X} s^{(\ell)}_j \cdot g_j^{(\ell)}(x)\right)$
            \State $\mathcal{S} \gets \mathcal{S} \cup \{(\ell,j)\}$
        \EndIf
    \EndIf
  \EndFor
  \State $b^{(\ell)} \gets b^{(\ell)} + \Delta^{(\ell)}$ \Comment{Apply bias shift}
  \State $\Delta h^{(\ell)} = \frac{1}{n_X} \sum_{k=1}^{n_X} \left( \relu\left(g^{(\ell)}(x^{(k)}) + \Delta^{(\ell)}\right) - \relu\left(g^{(\ell)}(x^{(k)})\right)  \right)$ 
  \State $b^{(\ell+1)} \gets b^{(\ell+1)} -  W^{(\ell+1)} \Delta h^{(\ell)}$ \Comment{Error correction}
  \State Recompute $g^{(\ell+1)}(x^{(k)})$, $h^{(\ell+1)}(x^{(k)}), L^{(\ell+1)}, U^{(\ell+1)}$ \Comment{Update}
\EndFor
\If{$d = \text{off}$}
    \For{$(\ell,j) \in \mathcal{S}$}
        \State Remove neuron $j$ from layer $\ell$
    \EndFor
\EndIf
\State Return $\{W^{(\ell)}, b^{(\ell)}\}_{\ell=1}^{L}, \ \{L^{(\ell)}, U^{(\ell)}\}_{\ell=1}^{L}, \ \mathcal{S}$ 
\end{algorithmic}
\end{algorithm}

\subsubsection{Pattern Merge}
\label{subsec:pm}

In the second step of the algorithm, neurons in the same layer with similar activation profiles are merged.

Related stability-based compression methods have already been proposed in the literature \citep{Kumar2019}. There, stably active neurons are merged with other stably active neurons based on linear dependencies in their weights. We extend this approach in two ways: 1) by merging neurons that are not necessarily stably active, and 2) by allowing neurons to be merged if their weights are only approximately linearly dependent. Specifically, two neurons are merged if the cosine similarity of their post-activation values exceeds a given threshold and their activation profile differs in less than a given percentage of samples in the probe batch. 

For all pairs $(i,j)$ of neurons that are not almost stable we calculate the cosine similarity of their centered post-activations and compare it to a threshold $\theta \le 1$, i.e.,
\begin{equation}\label{eq:merge_condition}
    \frac{\inner{\bar h^{(\ell)}_i(X),\,\bar h^{(\ell)}_j(X)}}{\norm{\bar h^{(\ell)}_i(X)}_{2}\,\norm{\bar h^{(\ell)}_j(X)}_{2}}\ge \theta,
\end{equation}
with $\bar h^{(\ell)}_i(X)$ being the centered post-activation of neuron $i$ obtained by subtracting the average post-activation, i.e.,
\begin{align}
    \bar h^{(\ell)}_i(X) := h_i^{(\ell)}(X) - \frac{1}{n_X} \inner{h_i^{(\ell)}(X), \mathbf{1}}.
\end{align}
This quantity is used because neuron outputs may exhibit different magnitudes or baseline activation levels while still being (approximately) linearly dependent.
If a pair of neurons $(i,j)$ fulfills \eqref{eq:merge_condition}, we treat it as a potential candidate for merging. If additionally their joint activation statistic exceeds a certain threshold $1 - \tau$, i.e., 
\begin{align}\label{eq:joint_activation_approximated}
a_{ij}^{(\ell)} := \frac{1}{n_X} \sum_{k=1}^{n_X} \left( \mathds{1}_{\{h_i^{(\ell)}(x^{(k)}) \geq 0\} \text{ and } \{h_j^{(\ell)}(x^{(k)}) \geq 0\}} \right) \geq 1 - \tau,
\end{align}
we merge the two neurons. For this, we calculate the coefficients 
\begin{align}\label{eq:alpha}
    \begin{split}
    \alpha_{ij}^{(\ell)} &=  \frac{\inner{h^{(\ell)}_i(X),\,h^{(\ell)}_j(X)}}{\norm{h^{(\ell)}_i(X)}_{2}^2}, \\
    \beta_{ij}^{(\ell)} &=  \frac{1}{n_X} \inner{h^{(\ell)}_j(X), \mathbf{1}} - \alpha_{ij}^{(\ell)} \frac{1}{n_X} \inner{ h^{(\ell)}_i(X), \mathbf{1}},
    \end{split}
\end{align}
to linearly express the output of neuron $j$ via a multiple of the output of neuron $i$ and an offset. The scaling factor $\alpha_j^{(\ell)}$ is optimal in the sense that it fulfills
\begin{equation}
    \alpha_{ij}^{(\ell)} \;\in\; \arg\min_{\alpha} \; \left\| h^{(\ell)}_j(X)-\alpha h^{(\ell)}_i(X) \right\|_{2}^2.
\end{equation}
Finally, neuron $j$ can be removed from layer $\ell$. Consequently, neuron $i$ now represents both contributions of both neurons $i$ and $j$. To account for this, we add a multiple of the $j$-th to the $i$-th column of $W^{(\ell+1)}$, i.e.,
\begin{equation}
    {W}^{(\ell+1)}_{:,i} \gets {W}_{:,i}^{(\ell+1)} + \alpha_{ij}^{(\ell)} {W}_{:,j}^{(\ell+1)}. 
\end{equation}
The approximation error of this operation can be estimated as
\begin{align}
    \Delta h^{(\ell)}_{ij} = \frac{1}{n_X} \sum_{k=1}^{n_X} \left( h_j^{(\ell)}(x^{(k)}) - \left(\alpha_{ij}^{(\ell)} h_i^{(\ell)}(x^{(k)}) + \beta_{ij}^{(\ell)} \right) \right),
\end{align}
which can be compensated by adapting the bias of the subsequent layer, i.e.,
\begin{align}
    \begin{split}
        b^{(\ell+1)} \gets b^{(\ell+1)} + \left(\beta_{ij}^{(\ell)} + \Delta h^{(\ell)}_{ij} \right)  {W}_{:,i}^{(\ell+1)}.
    \end{split}
\end{align}
The full pattern-merge procedure is described in \Cref{alg:patternmerge}.

\subsubsection{Complete algorithm}

\begin{figure}[H]
    \centering
    \begin{subfigure}{0.47\textwidth}
        \centering
        \includegraphics[width=\linewidth]{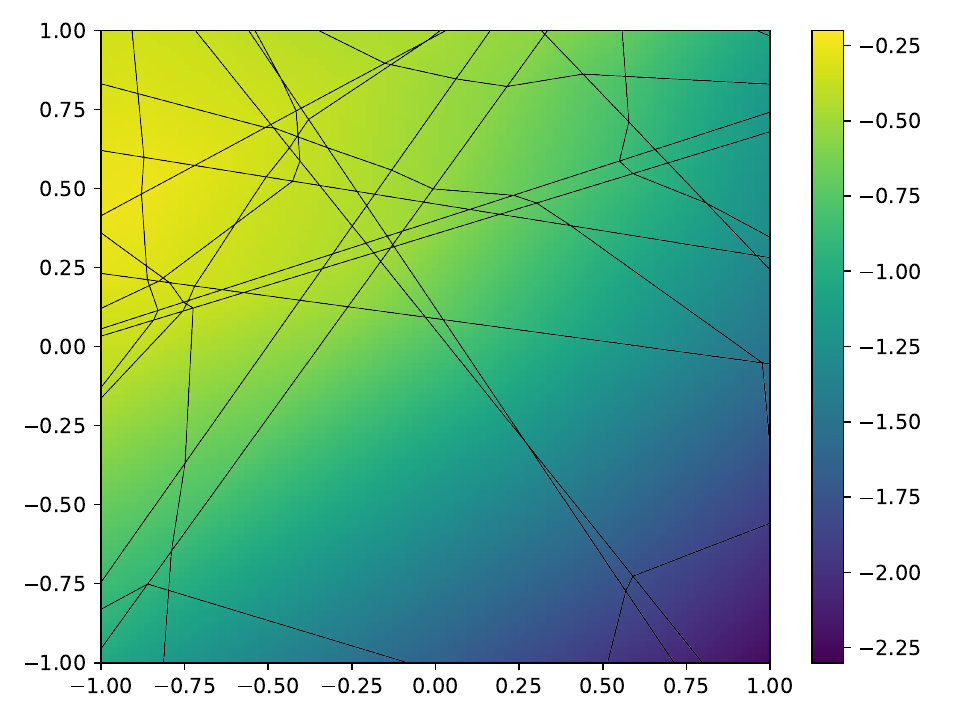}
        \caption{Original model (107 linear regions) with randomly initialized weights.}
        \label{fig:rere_0}
    \end{subfigure}
    \hfill
    \begin{subfigure}{0.47\textwidth}
        \centering
        \includegraphics[width=\linewidth]{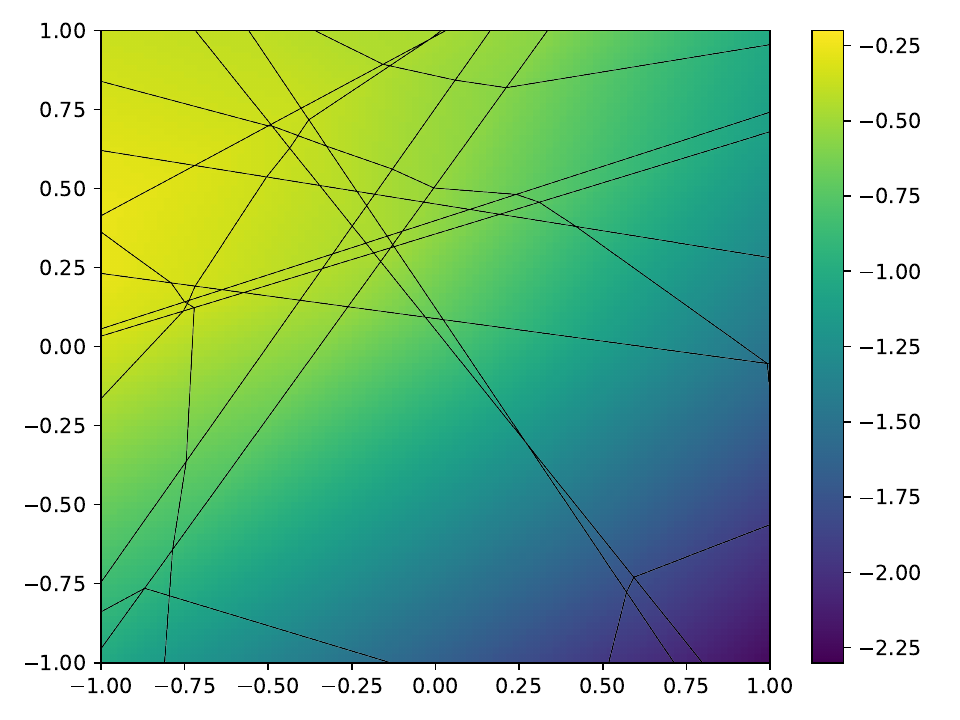}
        \caption{Intermediate model (75 linear regions) after stabilizing almost-off neurons via \textsc{BiasShift}. }
        \label{fig:rere_1}
    \end{subfigure}
    
    \vspace{1em}
    
    \begin{subfigure}{0.47\textwidth}
        \centering
        \includegraphics[width=\linewidth]{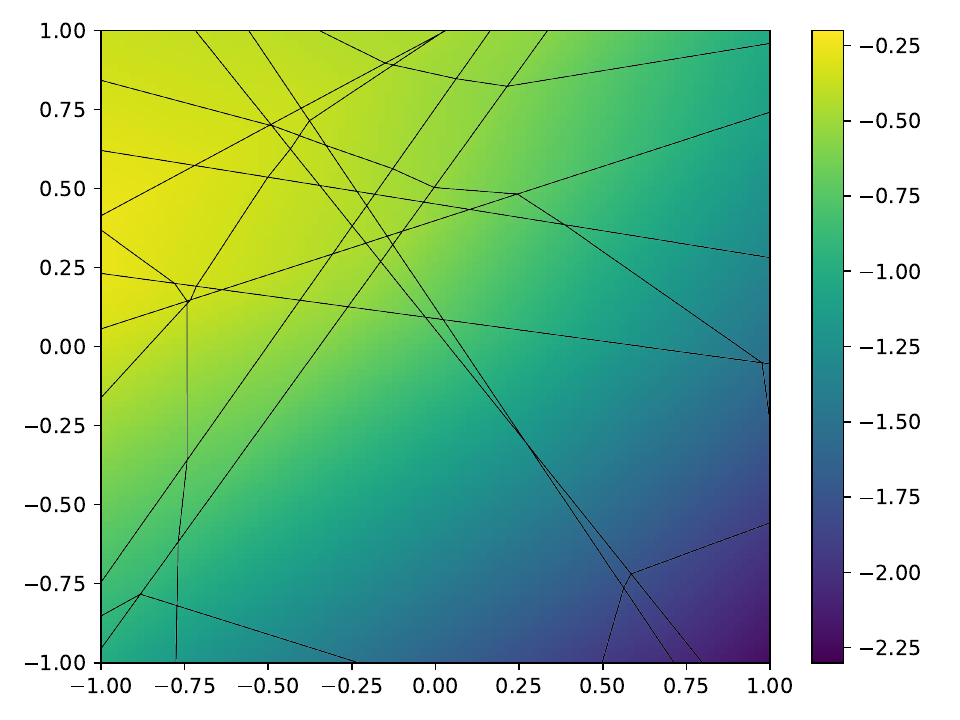}
        \caption{Intermediate model (65 linear regions) after merging similar neurons via \textsc{PatternMerge}.}
        \label{fig:rere_2}
    \end{subfigure}
    \hfill
    \begin{subfigure}{0.47\textwidth}
        \centering
        \includegraphics[width=\linewidth]{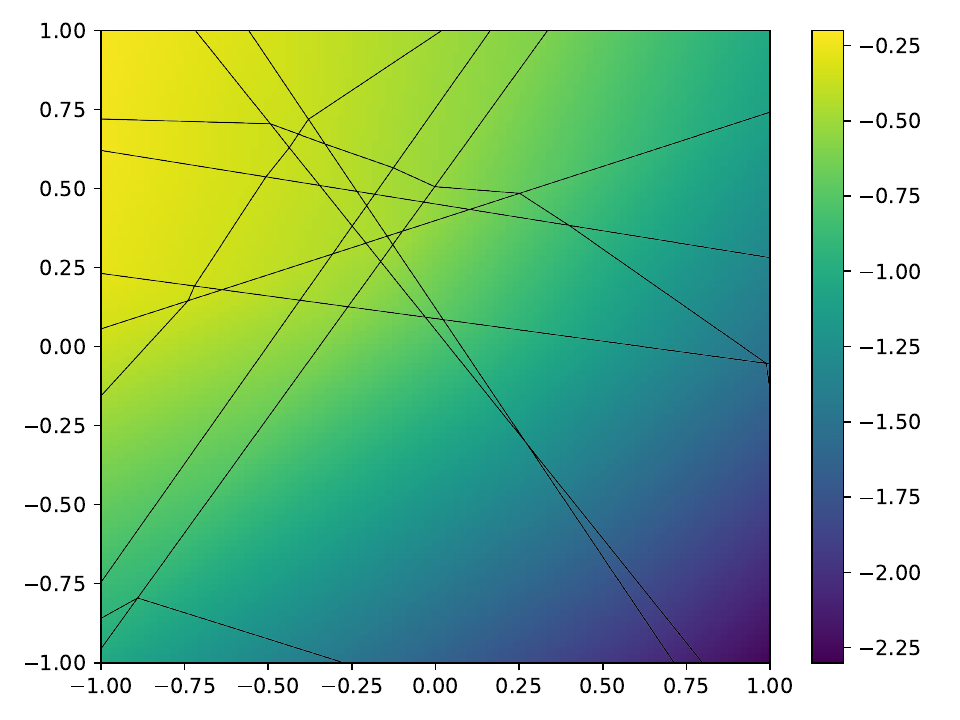}
        \caption{Region-reduced model (44 linear regions) after stabilizing almost-on neurons via \textsc{BiasShift}.}
        \label{fig:rere_3}
    \end{subfigure}
    \caption{Visualization of \textsc{RegionReduction} applied to a randomly initialized ANN with two hidden layers.}
    \label{fig:region_reduction}
\end{figure}

The complete algorithm is presented in \Cref{alg:region-reduction}. \Cref{fig:region_reduction} demonstrates the effect of its individual steps on a randomly initialized network with two hidden layers. It highlights the incremental simplification of the geometry of linear regions. The boundaries of the linear regions are highlighted as black lines and the background color reflects the output of the ANN. The effects of all substeps are evident, including the disappearance of decision boundaries due to the stabilization of almost stable neurons via \textsc{BiasShift}. In particular, \Cref{fig:rere_2} illustrates the collapse of two nearly parallel decision boundaries into one via \textsc{PatternMerge}. The region-reduced model retains a high level of agreement with the original model ($R^2 = 0.998$).

\begin{algorithm}[H]
\caption{\textsc{PatternMerge}}
\label{alg:patternmerge}
\begin{algorithmic}[1]
\Require Weights {$\{W^{(\ell)}, b^{(\ell)}\}_{\ell=1}^{L}$},  data {$X=\{x^{(k)}\}_{k=1}^{n_X}$}
\Require Merge threshold $\theta = 0.99$, activation safeguard $\tau = 0.05$.
\State Compute $\{h^{(\ell)}(X)\}_{\ell=1}^L$ \Comment{Initial forward pass}
\For{$\ell =1$ \textbf{to} $L-1$}
    \State $\mathcal P^{(\ell)} \gets \emptyset$ \Comment{Set of processed neurons}
    \ForAll{$(i,j) \in [n_\ell] \times [n_\ell]$ with $i < j$}
        \If{$i \in \mathcal P^{(\ell)}$ or $j \in \mathcal P^{(\ell)}$}
            \State continue
        \EndIf
        \State $\bar h_i^{(\ell)}(X) \gets h_i^{(\ell)}(X) - \frac{1}{n_X}  \inner{h^{(\ell)}_i(X), \ \mathbf{1}}$ \Comment{Centered post-activations}
        \State $\bar h_j^{(\ell)}(X) \gets  h_j^{(\ell)}(X) - \frac{1}{n_X}  \inner{h^{(\ell)}_j(X), \mathbf{1}}$ 
        \State $C_{ij} \gets \frac{\inner{\bar h^{(\ell)}_i(X),\,\bar h^{(\ell)}_j(X)}}{\norm{\bar h^{(\ell)}_i(X)}_{2}\,\norm{\bar h^{(\ell)}_j(X)}_{2}} $ \Comment{Cosine similarity}
        \If{$C_{ij}\ge \theta$} \Comment{Cosine similarity threshold}
                \If{$\frac{1}{n_X} \sum_{k=1}^{n_X}  \mathds{1}_{\{h_i^{(\ell)}(x^{(k)}) \geq 0 \text{ and } h_j^{(\ell)}(x^{(k)}) \geq 0\}} \geq 1 - \tau$} 
                    \State $ \mathcal P^{(\ell)} \gets \mathcal P^{(\ell)} \cup \{j\}$
                    \State $\alpha_{ij}^{(\ell)} \;=\;  \frac{\inner{h^{(\ell)}_i(X),\,h^{(\ell)}_j(X)}}{\norm{h^{(\ell)}_i(X)}_{2}^2}$
                    \State $\beta_{ij}^{(\ell)} =  \frac{1}{n_X} \inner{h^{(\ell)}_j(X), \mathbf{1}} - \alpha_{ij}^{(\ell)} \frac{1}{n_X} \inner{ h^{(\ell)}_i(X), \mathbf{1}}$
                    \State ${W}^{(\ell+1)}_{:,i} \gets {W}_{:,i}^{(\ell+1)} + \alpha_{ij}^{(\ell)} {W}_{:,j}^{(\ell+1)}$ \Comment{Merge neurons}
                    \State $\Delta h^{(\ell)}_{ij} = \frac{1}{n_X} \sum_{k=1}^{n_X} \left( h_j^{(\ell)}(x^{(k)}) - \left(\alpha_{ij}^{(\ell)} h_i^{(\ell)}(x^{(k)}) + \beta_{ij}^{(\ell)} \right) \right)$ 
                    \State $b^{(\ell+1)} \gets b^{(\ell+1)} + \left(\beta_{ij}^{(\ell)} + \Delta h^{(\ell)}_{ij} \right)  {W}_{:,i}^{(\ell+1)}$ \Comment{Error correction}
                \EndIf
            \EndIf
    \EndFor
    \State Recompute $g^{(\ell+1)}(x^{(k)})$, $h^{(\ell+1)}(x^{(k)})$ \Comment{Update}
    \For{$i \in \mathcal P^{(\ell)}$}
        \State Remove neuron $i$ from layer $\ell$
    \EndFor
\EndFor
\State \Return $\{W^{(\ell)}, b^{(\ell)}\}_{\ell=1}^L$
\end{algorithmic}
\end{algorithm}

\begin{algorithm}[htb]
\caption{\textsc{RegionReduction}}
\label{alg:region-reduction}
\begin{algorithmic}[1]
\Require Weights {$\{W^{(\ell)}, b^{(\ell)}\}_{\ell=1}^{L}$},  bounds $L^{(0)}, U^{(0)}$, data {$X=\{x^{(k)}\}_{k=1}^{n_X}$}
\State Compute bounds $\{L^{(\ell)}, U^{(\ell)}\}_{\ell=1}^{L}$ via $\eqref{eq:interval_arithmetic}$ or \eqref{prob:obbt}
\State Remove neurons with $U^{(\ell)}_i < 0, \ \ell \in [L-1], \ i \in [n_\ell]$ from network
\State Stabilize and remove almost-off neurons from network via \textsc{BiasShift} with $d$ = off
\State Merge redundant neurons via \textsc{PatternMerge}
\State Recompute/update bounds $\{L^{(\ell)}, U^{(\ell)}\}_{\ell=1}^{L}$ 
\State Stabilize almost-on neurons  via \textsc{BiasShift} with $d$ = on and save set of activated neurons $\mathcal{S}$
\State \Return $\{W^{(\ell)}, b^{(\ell)}\}_{\ell=1}^{L}, \{L^{(\ell)}, U^{(\ell)}\}_{\ell=1}^{L}, \mathcal{S}$
\end{algorithmic}
\end{algorithm}

\subsubsection{Hyperparameter search}

\begin{figure}[htpb]
    \centering
    \includegraphics[width=.7\linewidth]{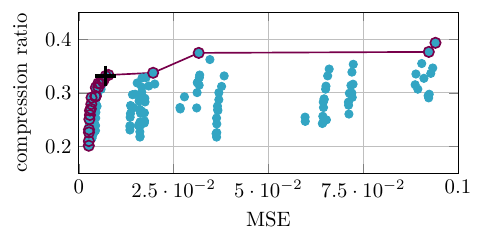}
    \caption{Pareto plot of averaged MSE and compression ratio for region-reduced ANNs, evaluated on ten representative networks from the four applications in \Cref{sec:applications}. The Pareto front is highlighted in red, and the selected hyperparameter configuration is marked with a black cross.} 
    \label{fig:hyperparameter}
\end{figure}

To calibrate the algorithm, we perform a hyperparameter search on ten exemplary networks from the four applications described in \Cref{sec:applications}. We evaluate the effect of the four parameters $\varepsilon, \rho, \theta, \tau$  on the region-reduced model's performance using a full-factorial study. Specifically, we compared the mean-squared error between the outputs of the region-reduced and original models, as well as the compression ratio of the region-reduced models. The hyperparameter values were drawn from the options given in \Cref{tab:hyperparameters}. Subsequently, Pareto-optimal hyperparameter combinations were identified, as shown in \Cref{fig:hyperparameter}. The selected hyperparameter configuration is used for all numerical examples in the following sections and can serve as a reasonable starting point for other applications. However, this is primarily a design choice for consistency, and a recalibration of the hyperparameters may be required for optimal performance in specific applications.

\begin{table}[t]
\centering
\caption{Investigated options in the hyperparameter study. The bold values are determined to be the best trade-off in terms of yielding models with low MSE and a high compression ratio.}
\label{tab:hyperparameters}
\small
\begin{tabular}{l p{4cm} p{6cm}}
\toprule
Parameter & Values & Description \\
\midrule
$\varepsilon$ & $0.01,\ 0.03,\ 0.05,\ \textbf{0.10}$ & Threshold of fraction of samples to detect active/inactive neurons \\
\addlinespace
$\rho$ & $0.10,\ 0.20,\ \textbf{0.30}$  & Maximum fraction of pre-activation range allowed for bias shifts \\
\addlinespace
$\theta$ & $0.90,\ 0.95,\ \textbf{0.99},\ 0.995$ & Threshold for cosine similarity in pattern merge \\
\addlinespace
$\tau$ & $0.01,\ 0.03,\ \textbf{0.05},\ 0.10$ & Allowed fraction of samples with dissimilar activation  \\
\bottomrule
\end{tabular}
\end{table}

\subsection{Pruned Model} \label{subsec:pruned}

We apply structured magnitude pruning by removing the neurons with smallest absolute weight values on each layer. In contrast to region-reduced models with their variable compression ratios, we prune a fixed percentage of neurons on each layer. For the numerical study in \Cref{sec:exper} we used a value of 20\%.
The resulting pruned models often resulted in infeasible optimization problems and numerical difficulties. Therefore, we used the probe batch $X$ as inputs and the original model's predictions $h^{(L)}(x^{(k)})$ as targets for Adam \citep{Kingma2014} in a fine-tuning procedure. 
This fine-tuning was restricted to a maximum of 25 epochs in the interest of a reduced time of pre-processing and used the same regularization hyperparameter $\lambda$ as in the training of the original model.

\subsection{Student Model} \label{subsec:student}

Again in the interest of comparability, we apply knowledge distillation by training from scratch a ReLU ANN of the exact same size as the pruned model, i.e., reduced in size by 20\% compared to the corresponding original models. In contrast to the pruned model, the network is initialized randomly and trained until convergence with a maximum number of 500 epochs using Adam \citep{Kingma2014}. While a larger time for compressing the model can be expected, this approach provides a benchmark for the accuracy of the compressed models.

\section{Applications}\label{sec:applications}

In this section we describe the optimization problems of increasing complexity. 
In Section~\ref{sec_app_peaks} we start with the minimization of the two-dimensional peaks function as a well-known benchmark problem in nonlinear optimization. The main focus here lies on the possibility to count and visualize the linear regions.
In Section~\ref{sec_app_superstructure} we investigate a superstructure model for the production of sustainable aviation fuels (SAF) from the literature \citep{Klimek2026}. It involves three physico-chemical processes: reverse water-gas shift (RWGS), biomass gasification, and Fischer-Tropsch (FT) synthesis. We approach this case study in two steps.
First, we optimize directly over the individual ANN surrogate models to study optimization problems based on realistic physical processes.
Second, we solve the full superstructure optimization problem that embeds all three neural networks. This step focuses on a larger-scale problem involving additional non-network submodels and variables.

\subsection{Peaks function} \label{sec_app_peaks}

In the first study, the ground truth is given via the peaks test function defined as 
\begin{align}
    \begin{split}
        f(x) &= 3 (1-x_1)^2 e^{-x_1^2 - (x_2+1)^2}  \\
        & \quad - 10 ( \frac{x_1}{5} - x_1^3 - x_2^5 ) e^{-x_1^2 - x_2^2} - \frac{1}{3} e^{-(x_1+1)^2 - x_2^2}
    \end{split}
    \label{eq:peaks}
\end{align}
on the domain $\mathcal D=[-3,3]^2$.
We train neural networks with five hidden layers and 25 neurons per layer using a maximum of 500 iterations of Adam \citep{Kingma2014} minimizing the mean-squared error between the training data and the networks' predictions. The networks were trained both with and without $\ell_1$ regularization, with the regularization parameter chosen as $\lambda =10^{-5}$ in the regularized setting. Each training was repeated five times. The size of the probe batch $X$ used for the compression algorithms was chosen as $n_X = 20000$.
 
\subsection{Production of Sustainable Aviation Fuels} \label{sec_app_superstructure}

We then move our focus to physico-chemical processes and superstructure optimization with embedded ANNs \citep{Klimek2026}. Specifically, we describe chemical compounds and conversions between them, with the aim of synthesizing a target product. These conversions are typically assumed to occur at fixed operating conditions such as temperature and pressure and hence fixed yield, usually determined from literature or simulation, and constant during optimization. This limits the space in which to optimize, as these operating conditions are not necessarily optimal in the overall system. By embedding ANNs for several conversion processes, we explicitly include operating conditions in the optimization problem as decision variables. Previous analysis has shown that this method does indeed improve the objective beyond the traditional approach \citep{Klimek2026}. Our model comprises mass and energy balances, CO\textsubscript{2} balancing, heat integration, and economic calculations. Its objective is the minimization of total cost while adhering to limits on total emissions. Streams are modelled as mixtures, necessitating partial and total mass balances, which determine the model's mathematical structure as mixed-integer quadratically constrained programming problem (MIQCP).

The superstructure in which our ANNs are embedded determines optimal synthesis pathways for the production of SAF. SAF provides a way to defossilize aviation while maintaining the utilization of existing infrastructure. Our model encompasses various conversion processes from raw materials to FT kerosene. Biomass, captured CO\textsubscript{2}, and natural gas are included as carbon sources, and can be converted to syngas, i.e., a mixture of CO and H\textsubscript{2}, as the central intermediate. Biomass is subjected to gasification, where it is decomposed to CO, CO\textsubscript{2}, and H\textsubscript{2}, under the supply of steam and O\textsubscript{2}. CO\textsubscript{2} -- atmospheric or derived from point sources -- is converted to CO via RWGS. Alternatively, syngas can be produced from natural gas via autothermal reforming. Acid gas removal is required to remove contaminants from raw syngas; electrolysis and air separation provide H\textsubscript{2} and O\textsubscript{2}. Syngas is then fed to a FT synthesis process, which yields n-alkanes distributed across an array of chain lengths. The fraction corresponding to kerosene is defined as the target compound. Carbon sequestration can remove CO\textsubscript{2} from the system and thereby provide offsets for any emissions, including supply chain emissions associated with feedstocks. The objective includes capital expenditure as well as the cost of raw materials, electricity, and heat. The complete model, all equations, and the corresponding data set are detailed in \citet{Klimek2026}. Embedded ANNs excluded, the problem comprises 5,325 continuous and 45 binary variables as well as 2,656 constraints. 

Three processes which exhibit significant variation depending on operating conditions are embedded as ANNs, namely RWGS, gasification, and FT~\citep{Klimek2026}, enabling us to capture the processes' complexity and dependence on operating conditions. RWGS converts CO\textsubscript{2} to CO using H\textsubscript{2}; the yield and energy requirements depend on the reactor temperature and the ratio of H\textsubscript{2} to CO\textsubscript{2} in the feed. The product of gasification varies significantly depending on the oxidizing agent, specifically the amount of steam (H\textsubscript{2}O) and O\textsubscript{2} fed to the reactor, as well as the supplemental CO\textsubscript{2}. These degrees of freedom can shift reactor operation between CO, CO\textsubscript{2}, and H\textsubscript{2} as products. FT yields a distribution of hydrocarbons, comprising light gases, gasoline, kerosene, diesel, and waxes. The kerosene fraction forms the target of the optimization problem. The distribution strongly depends on the reactor temperature and pressure as well as the composition of the syngas feed. Higher temperature, lower pressure and a higher fraction of H\textsubscript{2} shift it towards light gases, while lower temperature, higher pressure, and lower H\textsubscript{2} mass fraction favor the diesel fraction and waxes. Embedding the ANNs enables a representation of these dependencies in the superstructure, and simultaneous optimization of the topology, stream compositions, and operating conditions. 

Our initial investigation focuses on the maximization of the outputs of individual ANNs with specific objective functions for each process, i.e., we maximize the outlet CO mass fraction for RWGS and gasification, and the sum of kerosene-range mass fractions for FT. 

\begin{table}[!ht]
  \centering
  \small
  \caption{Input dimensions $n_0$ and output dimensions $n_L$ of the individual ANNs representing the three considered chemical processes, the corresponding sizes of the probe batch $X$ and the $\ell_1$ regularization parameters $\lambda$ used in this study.}
    \begin{tabular}{p{2cm}r p{4.8cm} ccc} 
    \toprule
      Process & $n_{0}$ & inputs & $n_L$ & $n_X$ & $\lambda$ \\%
    \midrule
    \multirow{2}{*}{RWGS}  & \multirow{2}{*}{2}     & reactor temperature & \multirow{2}{*}{8}     &  \multirow{2}{*}{20,000} & \multirow{2}{*}{$\{0, 10^{-4}\}$}\\
                       &       & H\textsubscript{2} mass fraction in &       &  \\
    \midrule
    \multirow{5}{*}{Gasification}    & \multirow{5}{*}{7}     & steam-to-biomass ratio                & \multirow{5}{*}{17}    & \multirow{5}{*}{200,000} & \multirow{5}{*}{$\{0, 10^{-5} \}$}\\%
               &       & CO\textsubscript{2}-to-biomass ratio  &       &  \\
                    &       & O\textsubscript{2}-to-biomass ratio   &       & \\
                    &       & gasifier temperature                  &       & \\
                    &       & biomass type (1 of 3 options)         &       & \\
    \midrule
    \multirow{3}{*}{FT}   & \multirow{3}{*}{3} & reactor temperature   & \multirow{3}{*}{38}    &  \multirow{3}{*}{50,000} & \multirow{3}{*}{$\{0, 10^{-4}\}$}\\
                            &   & reactor pressure      &       &  \\
                                &   & H\textsubscript{2} mass fraction in &       &  \\
    \bottomrule
    \end{tabular}%
  \label{tab:ANN_dimensions}%
\end{table}%

The data used for the training of the original ANNs is sampled from ASPEN Plus simulations \citep{Aspen1979}, with the number of data points ranging from 10,000 to 225,000. The input and output dimensions of the ANNs and the probe batch sizes are summarized in \autoref{tab:ANN_dimensions}. We use different-sized networks for the two kinds of optimization problems. For the optimization of individual ANNs we use larger networks than in the superstructure problem due to its increased complexity. The network sizes are illustrated in \Cref{tab:sizes_ANN_chemical}.

\begin{table}[htbp]
    \centering
    \small
    \caption{Sizes of ANNs embedded in the superstructure optimization problem and used for optimization of standalone networks.}
    \begin{tabular}{ccccc}
        \toprule
        \multirow{2}{*}{Optimization setting} &  \multirow{2}{*}{\# layers} & \multicolumn{3}{c}{\# neurons per layer} \\
        &         & FT & RWGS & Gasification \\
        \midrule
        \multirow{3}{*}{Individual Processes} & 4 & 30 & 30 & 30 \\
        & 6 & 30 & 30 & 30 \\
        & 8 & 30 & 30 & 30 \\
        \midrule
        \multirow{2}{*}{Superstructure} & 1 &    100 & 25 & 75 \\
        & 2 &    50 & 13 & 38 \\
        \bottomrule
    \end{tabular}
    \label{tab:sizes_ANN_chemical}
\end{table}
In the SAF superstructure problem, the ReLU ANNs with one hidden layer add a total of 688 continuous and 200 binary variables as well as 1,372 constraints to the problem, thereby multiplying the number of binary variables by four and increasing the number of constraints by 50\%.

\section{Numerical results and discussion}\label{sec:exper}

\begin{figure}[t]
\centering
    \includegraphics[width=0.49\textwidth]{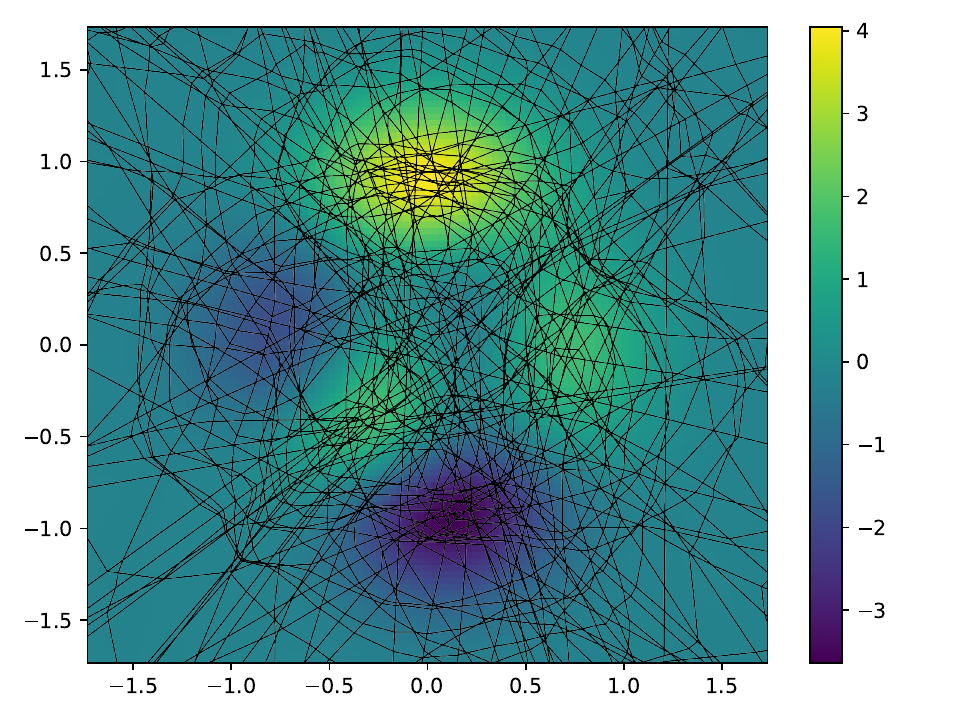}
\hfill
    \includegraphics[width=0.49\textwidth]{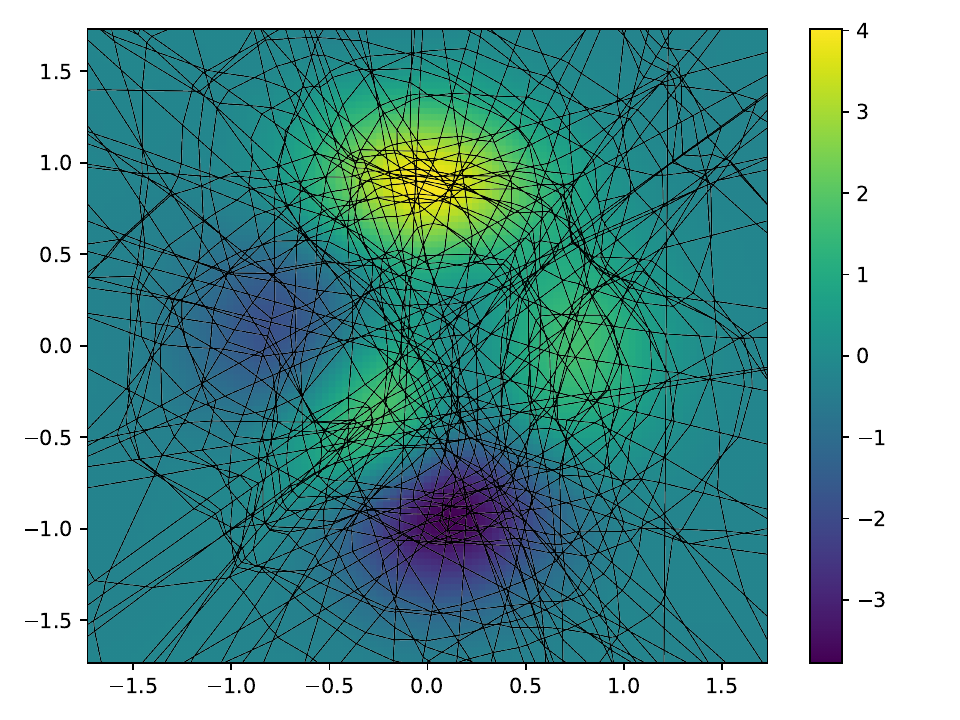}
    \includegraphics[width=0.49\textwidth]{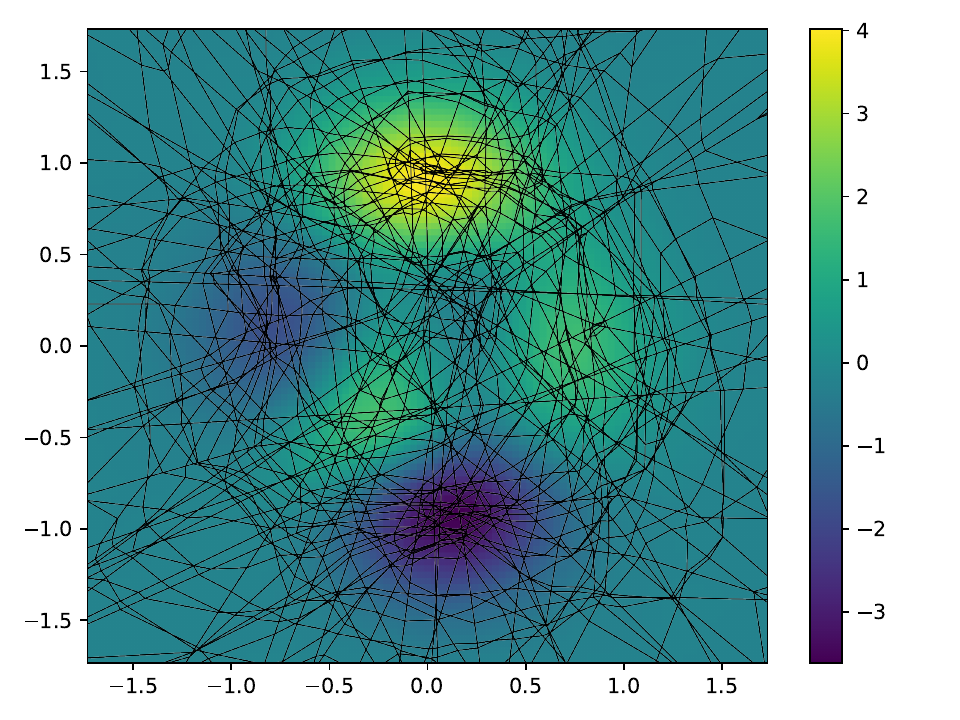}
\hfill
    \includegraphics[width=0.49\textwidth]{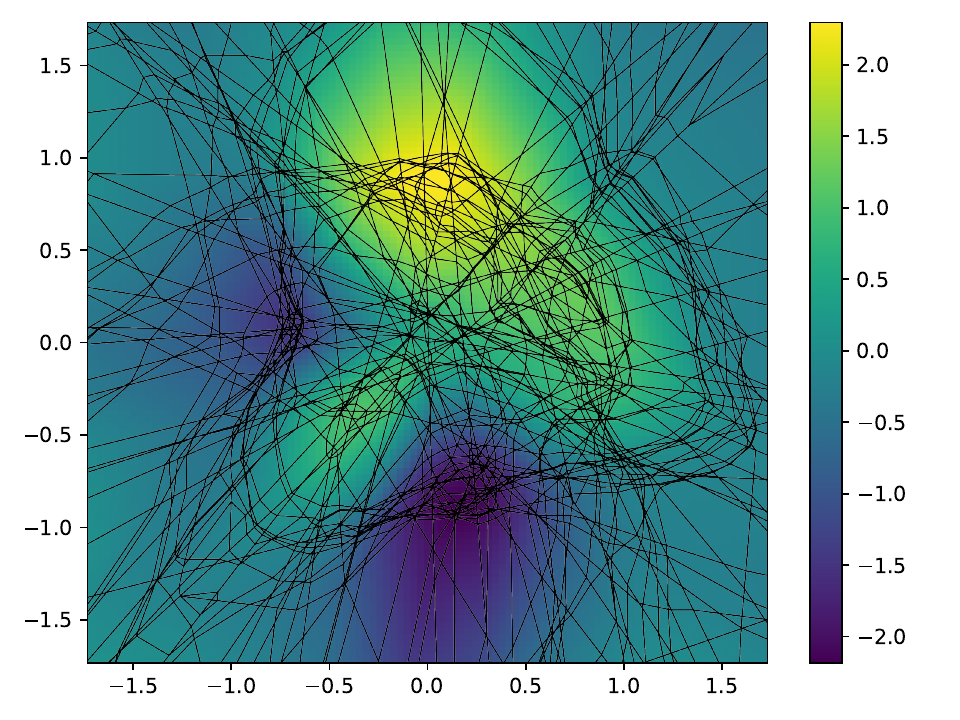}
\caption{Visualization of the linear regions for four models over the (rescaled) optimization domain $\mathcal{D}$. The color encoding shows resulting output values. The original model (top left, 5718 linear regions) is approximated by the three compressed models (region-reduced top right with 4716 linear regions, student bottom left with 4877 linear regions, pruned bottom right with 5114 linear regions). The region-reduced and student model manage to maintain the distinctive features of the peaks function, e.g., the positions and magnitudes of extrema with a visibly reduced number of linear regions. Magnitude pruning leads to a loss in accuracy.}
\label{fig:peaks_bspm}
\end{figure}

The procedure for obtaining numerical results is identical for all applications. First, the original ReLU ANN surrogate is trained using data from a ground truth model. The neural networks are implemented using Keras \citep{chollet2015keras} and trained with the optimization algorithm Adam \citep{Kingma2014}. Then, the corresponding region-reduced, pruned, and student models are derived. Lastly, the neural networks are translated into their MILP representation following \eqref{eq:bigM} using OMLT \citep{Ceccon2022}, and the resulting mixed-integer programs are solved via Gurobi v12.0 \citep{gurobi}. 

In all case studies we shall compare the four different models introduced in \Cref{sec:method} with respect to their predictive accuracy, the result of the optimization, and the overall computational time. The novel \Cref{alg:region-reduction} is used with bounds $\{L^{(\ell)}, U^{(\ell)}\}_{\ell=1}^{L}$ obtained from solving the bound-tightening LPs \eqref{prob:obbt}, as it allows a safe removal of more proven inactive neurons compared to using \eqref{eq:interval_arithmetic}. When formulating the MILP (or MIQCP in the superstructure problem), there are different ways to proceed. The first option is to disregard the bounds computed in \Cref{alg:region-reduction} and rely entirely on OMLT and Gurobi to handle the computation of bounds. The second option is to incorporate the pre-computed bounds $\{L^{(\ell)}, U^{(\ell)}\}_{\ell=1}^{L}$ directly into the MILP/MIQCP formulation, as their inclusion incurs no additional computational cost. Another consideration is whether to enable Gurobi's internal bound-tightening capabilities via the \texttt{OBBT} parameter \citep{gurobi}, which can further refine variable bounds during presolve and throughout the optimization process. For the standalone optimization problems in \cref{sec_peaks} and \cref{sec_standalone} we use the bounds from \Cref{alg:region-reduction} and leave the Gurobi parameter \texttt{OBBT} at its default value. For the superstructure problem in \cref{sec_superstructure} we rely on Gurobi's internal bound-tightening capabilities with the parameter \texttt{OBBT} set to 2.  

To verify that potential performance gains of the proposed algorithm in the applications in \cref{sec_peaks,sec_standalone} are not solely due to tighter bounds, we solve additional MILPs. More specifically, with region-reduced models the optimization problems are solved twice: once with setting the computed bounds $\{L^{(\ell)}, U^{(\ell)}\}_{\ell=1}^{L}$ when formulating the MILP and once without. For a meaningful comparison, the MILPs with the remaining models are also solved twice, once with default bounds and once with tighter bounds obtained from solving \eqref{prob:obbt}. This comparison is essential, as OBBT represents the natural first step to increase the computational efficiency of MILP with embedded ReLU ANNs. 

\subsection{Peaks function} \label{sec_peaks}

\begin{figure}[!ht]
    \centering
    \includegraphics[width=.95\linewidth]{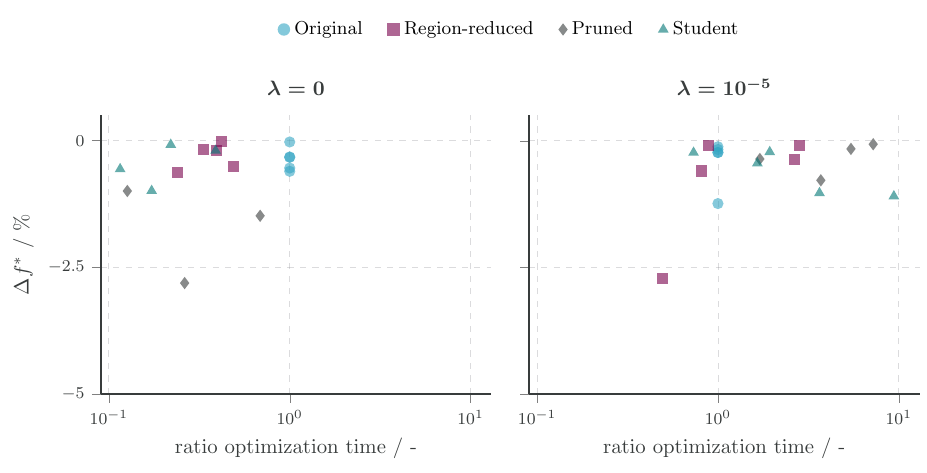}
    \caption{Deviations of objective values over ratios of $T_{\textrm{opt}}$ of compressed and original models for Peaks function.}
    \label{fig:peaks_table_illustrated}
\end{figure}

As the network input dimension is only 2, we apply a tool \citep{Plate2026} to visualize and count the numbers of linear regions. 
\Cref{fig:peaks_bspm} highlights the differences between the four considered models with respect to activation patterns for one of the five instances without regularization. The number of kinks in these black lines is bounded by the number of layers, in this case 5. In each region, the activation pattern is fixed, i.e., the model output varies linearly in the input $x$. We see a clear reduction in the number of linear regions of the region-reduced model compared to the original model, whereas the pruned model is not able to maintain an accurate resolution of the original activation pattern or the original function. This observation is supported by the results in \Cref{tab:results_peaks}, which show that the region-reduced models consistently reduce the average number of linear regions by ca. 16\% and 18\% for the unregularized and regularized networks while largely preserving the predictive qualities. Moreover, we observe that the drop in the average $R^2$ on the test set with $R^2 < 0.99$ stems from single outliers in which the globally used hyperparameters of \Cref{alg:region-reduction} seem too permissive, as shown, e.g., by the outlier in \Cref{fig:peaks_table_illustrated} on the right-hand side.

\begin{table}[!ht]
\setlength{\tabcolsep}{2pt}  
\scriptsize
\centering
\caption{Averaged times $T_{\textrm{model}}$ for obtaining the compressed neural network from the original model in seconds, mean absolute error and $R^2$ on the test set, number of linear regions \#LR of the ANN, optimization times in seconds, relative error of optimal objective with respect to known global minimum $\Delta f^*$, relative deviation of optimal inputs with respect to known global optimum $\Delta x^*$, and the total time $T_\textrm{total}$ comprising $T_\textrm{model}$ and $T_{\textrm{opt}}$. Results of models marked with a * are obtained by including big-M coefficients from \eqref{prob:obbt} in the MILP.}
\label{tab:results_peaks}
\begin{tabular}{lll cccc rccr}
\toprule
ANN & $\lambda$ & Model  & $T_{\textrm{model}}$ & MAE & $R^2$ & \#LR  &  $T_{\textrm{opt}}$ & $\Delta f^*$ & $\Delta x^*$ & $T_{\textrm{total}}$\\
\midrule
5$\times$25 & 0 & Original & 0.0 & 0.0226 & 0.9997 & 6716.6 & 146.52 &  -0.0036 & 0.0302 & 146.52  \\
5$\times$25 & 0 & Original* & 1.2 & 0.0226 & 0.9997 & 6716.6  & 49.27 &  -0.0036 & 0.0302 & 50.45  \\
5$\times$25 & 0 &  Region-reduced & 2.7 & 0.1012 & 0.9875 & 5670.0 & 49.34 &  -0.0030 & 0.0271 & 52.02  \\
5$\times$25 & 0 &  Region-reduced*& 2.7 & 0.1012 & 0.9875 & 5670.0 & 20.28 &  -0.0030 & 0.0271 & 22.97  \\
5$\times$25 & 0 & Pruned & 91.5 & 0.4366 & 0.8642 & 5154.2 & 31.68 &  -0.2121 & 0.2305 & 123.22  \\
5$\times$25 & 0 & Pruned*& 92.3 & 0.4366 & 0.8642 & 5154.2  & 9.99  & -0.2121 & 0.2305 & 102.25  \\
5$\times$25 & 0 & Student & 764.8 & 0.0353 & 0.9994 & 4321.4 & 20.57 &  -0.0063 & 0.0387 & 785.40  \\
5$\times$25 & 0 & Student* & 765.5 & 0.0353 & 0.9994 & 4321.4  & 5.91  & -0.0063 & 0.0387 & 771.42  \\
\midrule
5$\times$25 & 1e-05 & Original & 0.0 & 0.0201 & 0.9998 & 8560.0 & 10.18 &  -0.0040 & 0.0322 & 10.18  \\
5$\times$25 & 1e-05 & Original* & 1.1 & 0.0201 & 0.9998 & 8560.0 & 3.33 &  -0.0040 & 0.0322 & 4.41  \\
5$\times$25 & 1e-05 &  Region-reduced & 2.5 & 0.0782 & 0.9859 & 7045.4 & 10.35 &  -0.0077 & 0.0395 & 12.83  \\
5$\times$25 & 1e-05 &  Region-reduced*& 2.5 & 0.0782 & 0.9859 & 7045.4 & 1.56 &  -0.0077 & 0.0395 & 4.05  \\
5$\times$25 & 1e-05 & Pruned & 110.5 & 0.4201 & 0.8748 & 4232.0 & 24.91 &  -0.0289 & 0.0649 & 135.45  \\
5$\times$25 & 1e-05 & Pruned* & 111.3 & 0.4201 & 0.8748 & 4232.0 & 8.22 & -0.0289 & 0.0649 & 119.47  \\
5$\times$25 & 1e-05 & Student& 1587.1 & 0.0312 & 0.9995 & 5422.0 & 20.04 &  -0.0060 & 0.0359 & 1607.19  \\
5$\times$25 & 1e-05 & Student* & 1587.9 & 0.0312 & 0.9995 & 5422.0  & 4.10 &  -0.0060 & 0.0359 & 1591.95  \\

\bottomrule
\end{tabular}
\end{table}

The results in \Cref{tab:results_peaks} also illustrate that the student model performs best in terms of accuracy among the compressed ANNs. In contrast, the accuracy of the pruned models is the lowest. The higher accuracy of the student model comes at the cost of increased computational time $T_{\textrm{model}}$ required for solving the distillation problem, i.e., the training of the student model. Similarly, fine-tuning of the pruned model adds computational cost, but the impact is limited by the small number of epochs.

In the unregularized case, all compressed models exhibit reduced optimization times $T_\textrm{opt}$ compared to the original model, with speedups of approximately factor three for the region-reduced model and five to seven for the pruned and the student models, respectively, as is evident in \Cref{fig:peaks_table_illustrated} on the left. Applying OBBT and using the tighter bounds in the MILP consistently reduces $T_\textrm{opt}$ by again a factor of approximately three. As expected, the optimization times in the regularized cases are generally lower. Here, the compressed models offer no observable speedup compared to the original models, except for the case when bounds from OBBT are used. With these tighter bounds used in the MILP, the region-reduced model has the lowest average optimization time.

The quality of found solutions is approximately the same for original, region-reduced and student models, i.e., these two compressed models yield solutions close to those of the original model and close to the known global optimum. In contrast, the coarseness of the pruned models leads to larger deviations from the correct optimum. 
Moreover, due to the relatively small effort of compressing the model via \Cref{alg:region-reduction}, the region-reduced model is the only compressed model that yields a net efficiency gain in the unregularized case, i.e., the total time $T_\mathrm{total}$, comprising compression and subsequent optimization, remains significantly lower compared to the optimization time with the original model.

\subsection{Individual Processes}\label{sec_standalone}

\Cref{tab:rwgs_results,tab:gasification_results} depict the different performance indicators of the optimization problems of the RWGS and gasification processes. For a more concise presentation, the results for the FT case have been moved to the \autoref{app:ft}. Unlike the peaks function, the exact minima of these optimization problems are not known and the physical ground-truth models are processes modeled in ASPEN Plus \citep{Aspen1979} and not available for the evaluation of solutions. Therefore, for evaluating the relative deviations $\Delta x^*$ and $\Delta f^*$, the corresponding original model with its obtained solution is taken as a reference, even if the optimization problem with the original model did not converge.  In all problems, we are maximizing outlet mass fractions of specific components. Therefore, a negative $\Delta f^*$ implies that the found solution yields a better objective value when evaluated with the original model and compared to the solution obtained with the original model itself. 
First of all, the results again confirm that compression of ReLU ANNs via region-reduction is by a large margin the fastest among the three compression methods considered in this study. As can be read from the accuracies, most compressed models retain high levels of agreement with the test data as indicated by the $R^2$ values being close to $1$. In contrast to the results for the peaks function in \Cref{tab:results_peaks}, the accuracies of the pruned models are comparable to the others.

\begin{table}[H]
\centering
\setlength{\tabcolsep}{2pt}  
\caption{Performance of standalone \textbf{RWGS} optimization using ANNs with varying layers (30 neurons per layer) and $\ell_1$ regularization. Columns follow \Cref{tab:results_peaks}, with solved instances and average MIP gap of unsolved instances replacing linear region count. Models marked with a * incorporate big-M coefficients from \eqref{prob:obbt}.
}
\scriptsize

\label{tab:rwgs_results}
\begin{tabular}{lll rcc crcccr}
\toprule
ANN & $\lambda$ & Model  & $T_{\textrm{model}}$ & MAE & $R^2$ & \# Solved  &  $T_{\textrm{opt}}$ & Gap & $\Delta f^*$ & $\Delta x^*$ & $T_{\textrm{total}}$  \\
\midrule
4$\times$30 & 0 & Original & 0.0 & 0.0004 & 1.0000 & 5 & 28.61 & 0.000 & 0.0000 & 0.0000 & 28.61  \\
4$\times$30 & 0 & Original* & 1.0 & 0.0004 & 1.0000 & 5 & 4.44 & 0.000 & 0.0000 & 0.0000 & 5.45  \\
4$\times$30 & 0 &  Region-reduced & 2.3 & 0.0023 & 0.9930 & 5 & 10.61 & 0.000 & 0.0002 & 0.0000 & 12.94  \\
4$\times$30 & 0 &  Region-reduced*& 2.3 & 0.0023 & 0.9930 & 5 & 1.51 & 0.000 & 0.0002 & 0.0000 & 3.84  \\
4$\times$30 & 0 & Pruned & 79.6 & 0.0066 & 0.9842 & 5 & 12.12 & 0.000 & 0.0089 & 0.0132 & 91.67  \\
4$\times$30 & 0 & Pruned* & 80.2 & 0.0066 & 0.9842 & 5 & 3.93 & 0.000 & 0.0089 & 0.0132 & 84.12  \\
4$\times$30 & 0 & Student & 1273.0 & 0.0004 & 1.0000 & 5 & 1.80 & 0.000 & 0.0005 & 0.0000 & 1274.81  \\
4$\times$30 & 0 & Student* & 1273.6 & 0.0004 & 1.0000 & 5 & 0.65 & 0.000 & 0.0005 & 0.0000 & 1274.29  \\
\midrule
4$\times$30 & 0.0001 & Original & 0.0 & 0.0007 & 0.9998 & 5 & 0.09 & 0.000 & 0.0000 & 0.0000 & 0.09  \\
4$\times$30 & 0.0001 & Original* & 0.5 & 0.0007 & 0.9998 & 5 & 0.09 & 0.000 & 0.0000 & 0.0000 & 0.57  \\
4$\times$30 & 0.0001 &  Region-reduced & 1.4 & 0.0026 & 0.9937 & 5 & 0.05 & 0.000 & 0.0012 & 0.0010 & 1.43  \\
4$\times$30 & 0.0001 &  Region-reduced*& 1.4 & 0.0026 & 0.9937 & 5 & 0.05 & 0.000 & 0.0028 & 0.0039 & 1.43  \\
4$\times$30 & 0.0001 & Pruned & 97.1 & 0.0061 & 0.9839 & 5 & 9.83 & 0.000 & 0.0095 & 0.0167 & 106.92  \\
4$\times$30 & 0.0001 & Pruned* & 97.7 & 0.0061 & 0.9839 & 5 & 2.64 & 0.000 & 0.0095 & 0.0167 & 100.38  \\
4$\times$30 & 0.0001 & Student & 2253.5 & 0.0009 & 0.9998 & 5 & 0.06 & 0.000 & 0.0006 & 0.0019 & 2253.61  \\
4$\times$30 & 0.0001 & Student* & 2253.9 & 0.0009 & 0.9998 & 5 & 0.05 & 0.000 & 0.0006 & 0.0019 & 2253.97  \\
\midrule
6$\times$30 & 0 & Original & 0.0 & 0.0004 & 1.0000 & 5 & 728.61 & 0.000 & 0.0000 & 0.0000 & 728.61  \\
6$\times$30 & 0 & Original* & 3.9 & 0.0004 & 1.0000 & 5 & 137.20 & 0.000 & 0.0000 & 0.0000 & 141.11  \\
6$\times$30 & 0 & Region-reduced & 7.6 & 0.0051 & 0.9854 & 5 & 75.65 & 0.000 & 0.0002 & 0.0000 & 83.22  \\
6$\times$30 & 0 & Region-reduced*& 7.6 & 0.0051 & 0.9854 & 5 & 27.96 & 0.000 & 0.0002 & 0.0000 & 35.54  \\
6$\times$30 & 0 & Pruned & 103.8 & 0.0049 & 0.9940 & 5 & 174.18 & 0.000 & 0.0026 & 0.0000 & 277.95  \\
6$\times$30 & 0 & Pruned* & 106.1 & 0.0049 & 0.9940 & 5 & 49.52 & 0.000 & 0.0026 & 0.0000 & 155.63  \\
6$\times$30 & 0 & Student & 1012.0 & 0.0005 & 0.9999 & 5 & 57.18 & 0.000 & 0.0004 & 0.0000 & 1069.22  \\
6$\times$30 & 0 & Student* & 1014.2 & 0.0005 & 0.9999 & 5 & 12.68 & 0.000 & 0.0004 & 0.0000 & 1026.84  \\
\midrule
6$\times$30 & 0.0001 & Original & 0.0 & 0.0008 & 0.9998 & 5 & 0.24 & 0.000 & 0.0000 & 0.0000 & 0.24  \\
6$\times$30 & 0.0001 & Original* & 1.3 & 0.0008 & 0.9998 & 5 & 0.11 & 0.000 & -0.0000 & 0.0000 & 1.36  \\
6$\times$30 & 0.0001 & Region-reduced & 2.7 & 0.0034 & 0.9898 & 5 & 0.09 & 0.000 & 0.0005 & 0.0001 & 2.78  \\
6$\times$30 & 0.0001 & Region-reduced*& 2.7 & 0.0034 & 0.9898 & 5 & 0.08 & 0.000 & 0.0028 & 0.0009 & 2.78  \\
6$\times$30 & 0.0001 & Pruned & 130.6 & 0.0053 & 0.9928 & 5 & 36.03 & 0.000 & 0.0051 & 0.0060 & 166.59  \\
6$\times$30 & 0.0001 & Pruned* & 132.6 & 0.0053 & 0.9928 & 5 & 4.77 & 0.000 & 0.0051 & 0.0060 & 137.41  \\
6$\times$30 & 0.0001 & Student & 2947.1 & 0.0010 & 0.9997 & 5 & 0.11 & 0.000 & 0.0004 & 0.0001 & 2947.20  \\
6$\times$30 & 0.0001 & Student* & 2948.0 & 0.0010 & 0.9997 & 5 & 0.08 & 0.000 & 0.0004 & 0.0001 & 2948.05  \\
\midrule
8$\times$30 & 0 & Original & 0.0 & 0.0005 & 1.0000 & 1 & 6716.33 & 0.438 & 0.0000 & 0.0000 & 6716.33  \\
8$\times$30 & 0 & Original* & 10.7 & 0.0005 & 1.0000 & 2 & 3405.48 & 0.125 & -0.0000 & 0.0000 & 3416.21  \\
8$\times$30 & 0 & Region-reduced & 21.4 & 0.0070 & 0.9862 & 5 & 1864.49 & 0.000 & 0.0007 & 0.0000 & 1885.85  \\
8$\times$30 & 0 & Region-reduced*& 21.4 & 0.0070 & 0.9862 & 5 & 721.53 & 0.000 & 0.0007 & 0.0000 & 742.88  \\
8$\times$30 & 0 & Pruned & 131.8 & 0.0038 & 0.9971 & 3 & 2314.25 & 0.114 & 0.0020 & 0.0005 & 2446.08  \\
8$\times$30 & 0 & Pruned* & 138.3 & 0.0038 & 0.9971 & 5 & 1601.96 & 0.000 & 0.0020 & 0.0005 & 1740.24  \\
8$\times$30 & 0 & Student & 978.6 & 0.0006 & 0.9999 & 5 & 1264.01 & 0.000 & 0.0011 & 0.0000 & 2242.61  \\
8$\times$30 & 0 & Student* & 984.8 & 0.0006 & 0.9999 & 5 & 255.52 & 0.000 & 0.0011 & 0.0000 & 1240.34  \\
\midrule
8$\times$30 & 0.0001 & Original & 0.0 & 0.0009 & 0.9998 & 5 & 0.39 & 0.000 & 0.0000 & 0.0000 & 0.39  \\
8$\times$30 & 0.0001 & Original* & 2.2 & 0.0009 & 0.9998 & 5 & 0.16 & 0.000 & 0.0000 & 0.0000 & 2.38  \\
8$\times$30 & 0.0001 & Region-reduced & 4.4 & 0.0048 & 0.9803 & 5 & 0.15 & 0.000 & 0.0002 & 0.0000 & 4.51  \\
8$\times$30 & 0.0001 & Region-reduced*& 4.4 & 0.0048 & 0.9803 & 5 & 0.10 & 0.000 & 0.0009 & 0.0020 & 4.46  \\
8$\times$30 & 0.0001 & Pruned & 160.8 & 0.0047 & 0.9956 & 5 & 115.82 & 0.000 & 0.0015 & 0.0000 & 276.60  \\
8$\times$30 & 0.0001 & Pruned* & 166.1 & 0.0047 & 0.9956 & 5 & 19.75 & 0.000 & 0.0015 & 0.0000 & 185.84  \\
8$\times$30 & 0.0001 & Student & 3581.0 & 0.0011 & 0.9997 & 5 & 0.32 & 0.000 & 0.0014 & 0.0000 & 3581.30  \\
8$\times$30 & 0.0001 & Student* & 3582.6 & 0.0011 & 0.9997 & 5 & 0.13 & 0.000 & 0.0014 & 0.0000 & 3582.71  \\
\bottomrule
\end{tabular}

\end{table}
\begin{table}[H]
\centering
\setlength{\tabcolsep}{2pt}  
\caption{Performance of standalone \textbf{gasification} optimization using ANNs with varying numbers of layers (30 neurons per layer) and levels of $\ell_1$ regularization $\lambda$. The columns follow \Cref{tab:rwgs_results}.}
\scriptsize

\label{tab:gasification_results}
\begin{tabular}{lll rcc crcccr}
\toprule
ANN & $\lambda$ & Model  & $T_{\textrm{model}}$ & MAE & $R^2$ & \# Solved  &  $T_{\textrm{opt}}$ & Gap & $\Delta f^*$ & $\Delta x^*$ & $T_{\textrm{total}}$  \\
\midrule
4$\times$30 & 0 & Original & 0.0 & 0.0008 & 0.9998 & 5 & 22.50 & 0.000 & 0.0000 & 0.0000 & 22.50  \\
4$\times$30 & 0 & Original* & 1.4 & 0.0008 & 0.9998 & 5 & 6.46 & 0.000 & 0.0000 & 0.0000 & 7.82  \\
4$\times$30 & 0 & Region-reduced & 3.1 & 0.0014 & 0.9994 & 5 & 1.77 & 0.000 & 0.0004 & 0.0138 & 4.88  \\
4$\times$30 & 0 & Region-reduced*& 3.1 & 0.0014 & 0.9994 & 5 & 1.03 & 0.000 & 0.0004 & 0.0138 & 4.14  \\
4$\times$30 & 0 & Pruned & 688.2 & 0.0051 & 0.9963 & 5 & 6.89 & 0.000 & 0.0095 & 0.0220 & 695.14  \\
4$\times$30 & 0 & Pruned* & 689.1 & 0.0051 & 0.9963 & 5 & 1.36 & 0.000 & 0.0095 & 0.0220 & 690.50  \\
4$\times$30 & 0 & Student & 13040.1 & 0.0024 & 0.9990 & 5 & 4.88 & 0.000 & 0.0061 & 0.0894 & 13044.98  \\
4$\times$30 & 0 & Student* & 13041.0 & 0.0024 & 0.9990 & 5 & 2.66 & 0.000 & 0.0061 & 0.0894 & 13043.67  \\
\midrule
4$\times$30 & 1e-05 & Original & 0.0 & 0.0009 & 0.9998 & 5 & 0.51 & 0.000 & 0.0000 & 0.0000 & 0.51  \\
4$\times$30 & 1e-05 & Original* & 1.1 & 0.0009 & 0.9998 & 5 & 0.37 & 0.000 & 0.0000 & 0.0000 & 1.43  \\
4$\times$30 & 1e-05 & Region-reduced & 3.1 & 0.0013 & 0.9994 & 5 & 0.19 & 0.000 & 0.0041 & 0.0092 & 3.34  \\
4$\times$30 & 1e-05 & Region-reduced*& 3.1 & 0.0013 & 0.9994 & 5 & 0.20 & 0.000 & 0.0041 & 0.0092 & 3.35  \\
4$\times$30 & 1e-05 & Pruned & 969.5 & 0.0025 & 0.9989 & 5 & 2.19 & 0.000 & 0.0000 & 0.0000 & 971.69  \\
4$\times$30 & 1e-05 & Pruned* & 970.4 & 0.0025 & 0.9989 & 5 & 0.61 & 0.000 & 0.0000 & 0.0000 & 971.00  \\
4$\times$30 & 1e-05 & Student & 20471.8 & 0.0012 & 0.9996 & 5 & 0.38 & 0.000 & 0.0000 & 0.0000 & 20472.18  \\
4$\times$30 & 1e-05 & Student* & 20472.6 & 0.0012 & 0.9996 & 5 & 0.28 & 0.000 & 0.0000 & 0.0000 & 20472.86  \\
\midrule
6$\times$30 & 0 & Original & 0.0 & 0.0009 & 0.9998 & 5 & 2300.55 & 0.000 & 0.0000 & 0.0000 & 2300.55  \\
6$\times$30 & 0 & Original* & 4.9 & 0.0009 & 0.9998 & 5 & 607.08 & 0.000 & 0.0000 & 0.0000 & 612.02  \\
6$\times$30 & 0 & Region-reduced & 7.5 & 0.0019 & 0.9994 & 5 & 37.12 & 0.000 & 0.0031 & 0.0551 & 44.65  \\
6$\times$30 & 0 & Region-reduced*& 7.5 & 0.0019 & 0.9994 & 5 & 12.58 & 0.000 & 0.0031 & 0.0551 & 20.11  \\
6$\times$30 & 0 & Pruned & 1026.7 & 0.0055 & 0.9960 & 5 & 106.72 & 0.000 & 0.0068 & 0.0697 & 1133.39  \\
6$\times$30 & 0 & Pruned* & 1029.4 & 0.0055 & 0.9960 & 5 & 28.01 & 0.000 & 0.0068 & 0.0697 & 1057.44  \\
6$\times$30 & 0 & Student & 13411.4 & 0.0026 & 0.9989 & 5 & 200.52 & 0.000 & 0.0078 & 0.0458 & 13611.92  \\
6$\times$30 & 0 & Student* & 13414.6 & 0.0026 & 0.9989 & 5 & 97.22 & 0.000 & 0.0078 & 0.0458 & 13511.81  \\
\midrule
6$\times$30 & 1e-05 & Original & 0.0 & 0.0011 & 0.9997 & 5 & 3.35 & 0.000 & 0.0000 & 0.0000 & 3.35  \\
6$\times$30 & 1e-05 & Original* & 3.0 & 0.0011 & 0.9997 & 5 & 1.06 & 0.000 & -0.0000 & 0.0000 & 4.07  \\
6$\times$30 & 1e-05 & Region-reduced & 7.4 & 0.0015 & 0.9994 & 5 & 0.90 & 0.000 & 0.0065 & 0.0325 & 8.34  \\
6$\times$30 & 1e-05 & Region-reduced*& 7.4 & 0.0015 & 0.9994 & 5 & 0.45 & 0.000 & 0.0065 & 0.0325 & 7.89  \\
6$\times$30 & 1e-05 & Pruned & 1236.8 & 0.0028 & 0.9987 & 5 & 28.64 & 0.000 & -0.0000 & 0.0000 & 1265.46  \\
6$\times$30 & 1e-05 & Pruned* & 1239.6 & 0.0028 & 0.9987 & 5 & 2.65 & 0.000 & -0.0000 & 0.0000 & 1242.28  \\
6$\times$30 & 1e-05 & Student & 27273.2 & 0.0014 & 0.9996 & 5 & 1.62 & 0.000 & 0.0007 & 0.0083 & 27274.82  \\
6$\times$30 & 1e-05 & Student* & 27275.2 & 0.0014 & 0.9996 & 5 & 0.70 & 0.000 & 0.0007 & 0.0083 & 27275.85  \\
\midrule
8$\times$30 & 0 & Original & 0.0 & 0.0011 & 0.9998 & 0 & 7200.03 & 2.100 & 0.0000 & 0.0000 & 7200.03  \\
8$\times$30 & 0 & Original* & 13.1 & 0.0011 & 0.9998 & 0 & 7200.04 & 0.839 & -0.0056 & 0.0003 & 7213.17  \\
8$\times$30 & 0 & Region-reduced & 19.2 & 0.0022 & 0.9992 & 5 & 134.63 & 0.000 & -0.0047 & 0.0595 & 153.82  \\
8$\times$30 & 0 & Region-reduced*& 19.2 & 0.0022 & 0.9992 & 5 & 34.26 & 0.000 & -0.0047 & 0.0595 & 53.44  \\
8$\times$30 & 0 & Pruned & 1245.3 & 0.0060 & 0.9951 & 5 & 2356.23 & 0.000 & -0.0028 & 0.0769 & 3601.51  \\
8$\times$30 & 0 & Pruned* & 1252.6 & 0.0060 & 0.9951 & 5 & 425.70 & 0.000 & -0.0028 & 0.0769 & 1678.26  \\
8$\times$30 & 0 & Student & 17278.4 & 0.0029 & 0.9989 & 3 & 4026.74 & 0.851 & -0.0020 & 0.1015 & 21305.18  \\
8$\times$30 & 0 & Student* & 17286.4 & 0.0029 & 0.9989 & 3 & 1628.39 & 0.394 & -0.0020 & 0.1015 & 18914.79  \\
\midrule
8$\times$30 & 1e-05 & Original & 0.0 & 0.0012 & 0.9997 & 5 & 19.19 & 0.000 & 0.0000 & 0.0000 & 19.19  \\
8$\times$30 & 1e-05 & Original* & 7.2 & 0.0012 & 0.9997 & 5 & 2.47 & 0.000 & -0.0000 & 0.0000 & 9.67  \\
8$\times$30 & 1e-05 & Region-reduced & 15.5 & 0.0017 & 0.9993 & 5 & 2.47 & 0.000 & 0.0056 & 0.0223 & 17.98  \\
8$\times$30 & 1e-05 & Region-reduced*& 15.5 & 0.0017 & 0.9993 & 5 & 0.80 & 0.000 & 0.0056 & 0.0223 & 16.30  \\
8$\times$30 & 1e-05 & Pruned & 1561.6 & 0.0033 & 0.9983 & 5 & 225.54 & 0.000 & 0.0003 & 0.0189 & 1787.11  \\
8$\times$30 & 1e-05 & Pruned* & 1569.5 & 0.0033 & 0.9983 & 5 & 39.64 & 0.000 & 0.0003 & 0.0189 & 1609.14  \\
8$\times$30 & 1e-05 & Student & 33993.2 & 0.0016 & 0.9995 & 5 & 12.08 & 0.000 & 0.0007 & 0.0380 & 34005.30  \\
8$\times$30 & 1e-05 & Student* & 33997.4 & 0.0016 & 0.9995 & 5 & 1.46 & 0.000 & 0.0007 & 0.0380 & 33998.84  \\
\bottomrule
\end{tabular}

\end{table}
Again, there are some outliers, e.g., the regularized instances with 8 layers of the RWGS process, where region-reduced models exhibit slightly worse $R^2$ values. As in the peaks example, these outliers can be traced back to individual ANNs.
Regarding the optimization times, there is again a clear divide between regularized and unregularized configurations. While optimization problems with the largest regularized networks are solvable within seconds, the same optimization problems with the unregularized original models are largely intractable, frequently hitting the time limit. For the unregularized configurations, all compressed models are able to reduce the optimization times. This speedup is most apparent with the largest models, where almost all optimization problems with compressed models are now tractable. Notably, the region-reduced compression is the only method examined in this study that enables successful optimization across all instances of the largest unregularized networks, regardless of whether bounds from LP-based bound-tightening are used. Across the unregularized models in both processes, the region-reduced model delivers net efficiency gains, whereas the student model does so only for the RWGS case in \Cref{tab:rwgs_results}, as the training in the gasification case is more expensive due to higher input/output dimensions and the larger training set.
For the regularized networks, speedups of $T_\mathrm{opt}$ can be observed for region-reduced and student models, whereas pruned models sometimes exhibit an increase in computational time, e.g., in  \Cref{tab:rwgs_results} for the regularized models with eight layers. 

The quality of obtained solutions with the compressed models is generally high, with all objective values showing less than 1\% deviation from the corresponding optimal objective of the original model. The largest measured average deviation $\Delta f^*$ in \Cref{tab:rwgs_results,tab:gasification_results} is 0.95\%, observed in regularized pruned models with 4 layers for the gasification problem. In the largest unregularized instance of the gasification optimization, the compressed models all yield negative $\Delta f^*$, indicating that the found solution is better than the solution obtained with the original model. This occurs as all original models fail to converge within the time limit, returning a suboptimal solution, whereas the compressed models are able to reach the global optimum.

\subsection{Superstructure Problem} \label{sec_superstructure}
We now consider the SAF superstructure problem, which comprises 14 instances of scenarios that vary in their restrictions on CO\textsubscript{2} emissions, among other parameters. As discussed before, the superstructure problem is formulated as a MIQCP, incorporating three embedded ReLU ANNs, along with other binary and continuous decision variables. We examine two architectural configurations for the ANNs: one in which each network contains a single hidden layer, and another in which each network comprises two hidden layers, see \Cref{tab:sizes_ANN_chemical}.  Furthermore, we investigate two regularization settings: one without $\ell_1$ regularization, and another in which $\ell_1$ regularization is applied to each ANN, with varying magnitudes $\lambda$ as described in \Cref{tab:ANN_dimensions}. 

\Cref{tab:results_superstructure} details the optimization results with the original and the compressed models. As the models are now used in multiple optimization problems, the computational effort of obtaining the compressed model in the first place is of less interest as before. Hence, we only report the average optimization time $T_{\textrm{opt}}$ in seconds. Also, as exact optimal solutions are again unknown, we define the configuration with the regularized 1-layer networks as the overall baseline for comparing the found solutions. This is done as this configuration is the one used in the original publication \citep{Klimek2026}. The results indicate that the region-reduced models decrease the average optimization time by approximately 40\% for regularized ANNs with a single hidden layer and by 50\% for those with two hidden layers compared to the respective original models. In contrast, pruning seems inadequate for this task, as the majority of instances cannot be solved within the time limit and those instances that are solved yield solutions with the largest deviations compared to the baseline solution. The student models yield moderate average speedups for unregularized networks, e.g., 40\% faster optimization for the single layer case, but offer no measurable advantage for the regularized single layer configuration.

\begin{figure}[htbp]
    \centering
    \includegraphics[width=.97\linewidth]{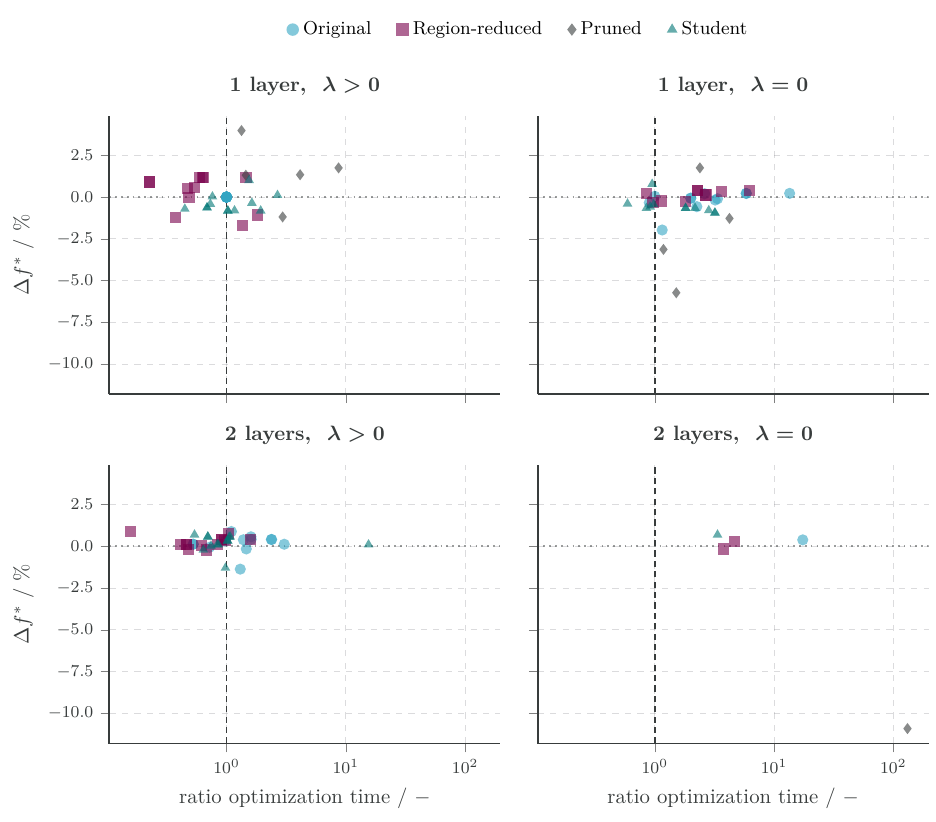}
    \caption{Relative deviation of optimal objective values over ratio of optimization times obtained with the compressed models compared to the baseline with regularized original models with one hidden layer for all instances of the superstructure problem considered.}
    \label{fig:superstructure_deltaf_1layerl1}
\end{figure}

Comparing the quality of found solutions as measured by the relative deviation of the optimal objectives to the baseline of the regularized one layer baseline, \Cref{fig:superstructure_deltaf_1layerl1} shows that most instances have a negligible relative deviation $\Delta f^*$. This is especially pronounced with the regularized networks with two layers, whereas the variance of $\Delta f^*$ is slightly higher for the smaller, regularized models. This suggests that the robustness to compression might increase with the number of layers. Also the slightly higher accuracies of the region-reduced networks with two layers hint in this direction.

\setlength{\tabcolsep}{2pt}  
\begin{table}[H]
  \centering
  \caption{Performance of ANNs  with varying numbers of layers and  $\ell_1$ regularization on the superstructure problem. Shown are mean test error, $R^2$, number of solved (S), timed-out (T), and infeasible (I) instances, average optimization time $T_{\textrm{opt}}$ on solved instances, geometric mean ratio (GMR) of optimization times relative to overall baseline (1-layer model with $\lambda$ > 0) and corresponding original model in each configuration. MIP gaps are averaged over unsolved instances.  Average relative deviations of the optimal objective value $\Delta f^*$ and the obtained minimizers $\Delta x^*$ compare the found solutions to the solution of the overall baseline. Averages indicated by a $^{*}$ are based on only one instance.} 
  \label{tab:results_superstructure}
\scriptsize
\begin{tabular}{@{}ll ll  cr cc lll}
\toprule
\multirow{2}{*}{Config.} & \multirow{2}{*}{Model} & \multirow{2}{*}{MAE} & \multirow{2}{*}{$R^2$} & Instances & \multirow{2}{*}{$T_{\textrm{opt}}$} & GMR  & GMR & \multirow{2}{*}{Gap} & \multirow{2}{*}{$\Delta f^*$}  & \multirow{2}{*}{$\Delta x^*$} \\
& & & &  S / T / I & & (all) & (config.)\\
\midrule
  & Original & 0.0005 & 0.9994 & 13 / 1 / 0  & 91.2 & (1.00) & (1.00)& 0.118$^{*}$ & 0.0000 & 0.0000\\
 1 layer & Region-reduced & 0.0012  & 0.9937 & 12 / 1 / 1  & 65.5 & 0.60 & 0.60 & 0.133$^{*}$ &  0.0030 & 0.0370 \\
 $\lambda > 0$ & Pruned & 0.0128 & 0.8832 & 5 / 8 / 1  & 319.6 & 2.90 & 2.90 & 0.128 & 0.0144 & 0.1045 \\
 & Student & 0.0006 & 0.9993 & 12 / 2 / 0  & 119.9 & 1.05 & 1.05 & 0.039 &-0.0040 & 0.0508 \\
\midrule
 & Original &0.0004  & 0.9997 & 12 / 2 / 0  & 265.8 & 2.66 & (1.00) & 0.306 & -0.0020 & 0.0484 \\
 1 layer & Region-reduced & 0.0012 & 0.9945 & 10 / 3 / 1  & 215.3 & 2.07 & 0.66 & 0.128 & 0.0012 & 0.0576\\
 $\lambda = 0$ & Pruned & 0.0113 & 0.9006 & 4 / 9 / 1  & 260.5  & 2.06 & 0.99 & 0.232 & -0.0210 & 0.0633 \\
 & Student & 0.0008 & 0.9995 & 11 / 3 / 0  & 158.1  & 1.50 & 0.58 & 0.197 & -0.0054 & 0.0736 \\
\midrule
& Original & 0.0006 &0.9996 & 13 / 1 / 0  & 113.4 & 1.28 & (1.00) & 0.158$^{*}$ &  0.0017 & 0.0367\\
2 layers & Region-reduced & 0.0009 & 0.9957 & 13 / 1 / 0  & 69.8 & 0.64 & 0.51 & 0.470$^{*}$ & 0.0025 & 0.0337 \\
$\lambda  > 0$ & Pruned & 0.0091 & 0.9304 &  0 / 13 / 1   & -- & -- & -- & 0.409 & -- & -- \\
 & Student & 0.0008 & 0.9995 & 13 / 1 / 0  & 181.1 & 1.04 & 0.81 & 0.197$^{*}$ & 0.0019 & 0.0422\\
\midrule
& Original & 0.0005  &0.9999 & 1 / 13 / 0  & 49.5$^{*}$ & 17.41$^{*}$ & (1.00)$^{*}$ & 0.586 & 0.0039$^{*}$ & 0.0226$^{*}$ \\
2 layers & Region-reduced & 0.0009 &0.9973 & 2 / 12 / 0  & 184.6 & 68.4 & 0.31 & 0.436 &  0.0005 & 0.0277 \\
$\lambda = 0$ & Pruned & 0.0124 & 0.9054 & 1 / 13 / 0  & 373.6$^{*}$ & 373.6$^{*}$ & 7.55$^{*}$ & 0.551 & -0.1090$^{*}$ & 0.0948$^{*}$ \\
& Student & 0.0012 & 0.9996 & 1 / 13 / 0 & 9.5$^{*}$ & 3.36$^{*}$ & 0.19$^{*}$ & 0.471 & 0.0069$^{*}$ & 0.0286$^{*}$ \\
\bottomrule
\end{tabular}
\end{table}

\begin{figure}[t]
    \centering
    \includegraphics[width=.97\linewidth]{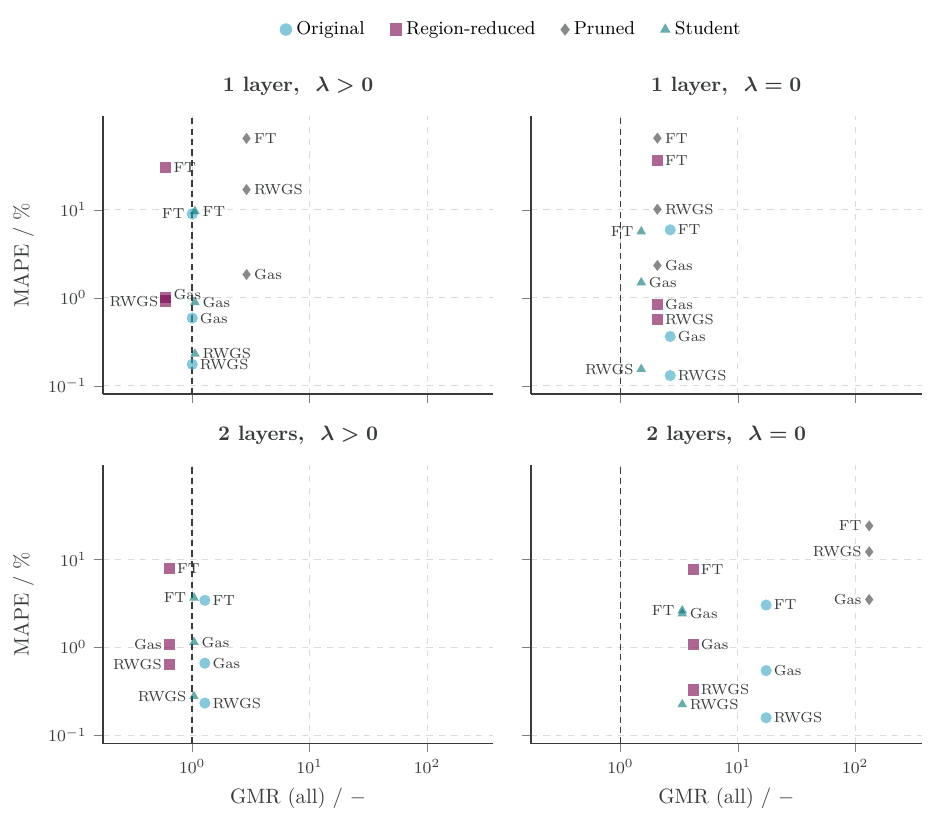}
    \caption{Mean-absolute percentage error over geometric mean of ratios of optimization times obtained with the compressed models compared to the baseline with regularized original models with one hidden layer for all solved instances of superstructure problems.}
    \label{fig:superstructure_mape_1layer_1}
\end{figure}

On the downside, embedding the region-reduced and pruned models leads to two and three infeasible instances, respectively. This behavior is most likely caused by the FT subprocess. As explained in \autoref{sec_app_superstructure}, the FT ANN predicts the mass fraction distribution across 30 hydrocarbon products of varying chain lengths. The superstructure problem uses a constraint that all predicted mass fractions sum to one. This constraint is enforced by computing the final mass fraction as the residual value required for the constraint to be satisfied. Consequently, small overestimations of the predicted mass fractions can easily lead to infeasibilities. This explanation is supported by the observation that the FT ANNs exhibit the highest prediction errors among the three processes, as illustrated in \Cref{fig:superstructure_mape_1layer_1}. 

Moreover, \Cref{fig:superstructure_mape_1layer_1} hints at another observation. Namely, that there seems to be a trade-off between optimization time and accuracy of the employed models. For slight losses in accuracy, gains in computational efficiency can be achieved, as seen with the region-reduced models. However, when model accuracy degrades further, this trend reverses and computational times may increase, as observed with the pruned models in the regularized configurations. This suggests that while some approximation error due to more compact models is tolerable or even beneficial for speed, excessive error complicates the satisfaction of constraints and ultimately degrades solver performance. \autoref{app:superstructure} contains the plots analogous to \Cref{fig:superstructure_deltaf_1layerl1,fig:superstructure_mape_1layer_1} which use each configuration's original model as the baseline.

\clearpage

\newpage

\section{Conclusion}

We presented a novel compression method for ReLU ANNs designed to increase the tractability of MILP formulations with embedded neural networks.
The two substeps of the algorithm, \textit{bias shift} and \textit{pattern merge}, stabilize previously unstable neurons and reduce activation redundancy while largely preserving functional fidelity, thereby lowering the number of its linear regions. Our method requires minimal computational effort. If OBBT is employed, the method's computational cost is dominated by the OBBT step, while avoiding the additional training required by the magnitude pruning and knowledge distillation approaches examined in this paper.
Across benchmarks, the region-reduced models consistently produce optimal solutions that closely match those of the original network or the ground-truth. Additionally, these solutions are obtained more quickly with speedups of 40\% to 50\% on the investigated superstructure problems, making the method a cost-effective pre-processing step for optimization problems with embedded ReLU neural networks. The results also indicate that both the modeling and compression of ReLU ANNs must be performed with care in the presence of additional constraints, such as those arising in the superstructure problem. 

The proposed method depends on the activation statistics of individual neurons, which are derived from calibration data sampled from the input domain $\mathcal{D}$. Consequently, the method's performance is influenced by both the sampling method and the sampling size. This may limit the applicability in high-dimensional settings and should be investigated in the future. Moreover, extending the framework to convolutional or hybrid architectures, and to alternative piecewise-linear activations (e.g., Leaky ReLU or Maxout), remains an open line of further research. 




\section*{CRediT authorship contribution statement}

\textbf{Christoph Plate}: Conceptualization, Methodology, Software, Formal Analysis, Investigation, Writing - original draft, Visualisation
\textbf{Caroline Ganzer}: Conceptualization, Supervision, Writing – review and editing
\textbf{Mirko Hahn}: Conceptualization
\textbf{Alexander Klimek}: Conceptualization, Formal analysis, Investigation, Writing – review and editing, Visualisation
\textbf{Heyuan Liu}: Conceptualization, Software, Investigation
\textbf{Sebastian Sager}: Resources, Funding acquisition, Writing – review and editing, Supervision
\textbf{Kai Sundmacher}: Resources, Funding acquisition, Supervision
\textbf{Hanna Wilhelm}: Conceptualization, Writing – review and editing

\section*{Declaration of competing interest}

The authors declare that they have no known competing financial interests or personal relationships that could have appeared to influence the work reported in this paper.

\section*{Data availability}

Data will be made available upon reasonable request.

\section*{Acknowledgements}

The work was funded by the German Research Foundation DFG within the priority program 2331 ’Machine Learning in Chemical Engineering’ under grant SA 2016/3-1 and by the European Regional Development Fund (ERDF) within the programme Research and Innovation - Grant Number ZS/2023/12/182075, which we gratefully acknowledge. We further acknowledge financial support from the research initiative ’SmartProSys: Intelligent Process Systems for the Sustainable Production of Chemicals’ funded by the Ministry for Science, Energy, Climate Protection and the Environment of the State of Saxony-Anhalt.

\bibliographystyle{elsarticle-harv}
\bibliography{reference,reference2}

\begin{appendices}

\section{FT results}\label{app:ft}

\Cref{tab:ft_results} illustrates the optimization results for the standalone FT optimization problem.

\begin{table}[htbp]
\centering
\setlength{\tabcolsep}{2pt}  
\scriptsize
\caption{Performance of standalone \textbf{FT} optimization using ANNs with varying numbers of layers (30 neurons per layer) and levels of $\ell_1$ regularization $\lambda$. The columns follow \Cref{tab:rwgs_results}.
}
\label{tab:ft_results}
\begin{tabular}{lll rcc crcccr}
\toprule
ANN & $\lambda$ & Model  & $T_{\textrm{model}}$ & MAE & $R^2$ & \#Solved  &  $T_{\textrm{opt}}$ & Gap & $\Delta f^*$ & $\Delta x^*$ & $T_{\textrm{total}}$  \\
\midrule
4$\times$30 & 0 & Original & 0.0 & 0.0004 & 0.9997 & 5 & 19.45 & 0.000 & 0.0000 & 0.0000 & 19.45  \\
4$\times$30 & 0 & Original* & 2.0 & 0.0004 & 0.9997 & 5 & 5.23 & 0.000 & 0.0000 & 0.0000 & 7.27  \\
4$\times$30 & 0 & Region-reduced & 4.7 & 0.0008 & 0.9938 & 5 & 5.58 & 0.000 & 0.0012 & 0.0276 & 10.33  \\
4$\times$30 & 0 & Region-reduced*& 4.7 & 0.0008 & 0.9938 & 5 & 1.76 & 0.000 & 0.0012 & 0.0276 & 6.50  \\
4$\times$30 & 0 & Pruned & 188.1 & 0.0018 & 0.9931 & 5 & 4.22 & 0.000 & 0.0082 & 0.1944 & 192.35  \\
4$\times$30 & 0 & Pruned* & 189.6 & 0.0018 & 0.9931 & 5 & 1.95 & 0.000 & 0.0082 & 0.1944 & 191.56  \\
4$\times$30 & 0 & Student & 2417.7 & 0.0004 & 0.9997 & 5 & 2.19 & 0.000 & 0.0037 & 0.0896 & 2419.89  \\
4$\times$30 & 0 & Student* & 2419.2 & 0.0004 & 0.9997 & 5 & 1.08 & 0.000 & 0.0037 & 0.0896 & 2420.26  \\
\midrule
4$\times$30 & 0.0001 & Original & 0.0 & 0.0007 & 0.9990 & 5 & 0.35 & 0.000 & 0.0000 & 0.0000 & 0.35  \\
4$\times$30 & 0.0001 & Original* & 0.9 & 0.0007 & 0.9990 & 5 & 0.39 & 0.000 & -0.0000 & 0.0000 & 1.29  \\
4$\times$30 & 0.0001 & Region-reduced & 2.2 & 0.0008 & 0.9972 & 5 & 0.17 & 0.000 & 0.0005 & 0.0141 & 2.36  \\
4$\times$30 & 0.0001 & Region-reduced*& 2.2 & 0.0008 & 0.9972 & 5 & 0.19 & 0.000 & 0.0004 & 0.0145 & 2.39  \\
4$\times$30 & 0.0001 & Pruned & 232.6 & 0.0018 & 0.9932 & 5 & 1.63 & 0.000 & 0.0076 & 0.1248 & 234.25  \\
4$\times$30 & 0.0001 & Pruned* & 233.9 & 0.0018 & 0.9932 & 5 & 0.99 & 0.000 & 0.0076 & 0.1248 & 234.93  \\
4$\times$30 & 0.0001 & Student & 5299.9 & 0.0007 & 0.9988 & 5 & 0.08 & 0.000 & 0.0054 & 0.0791 & 5300.00  \\
4$\times$30 & 0.0001 & Student* & 5300.4 & 0.0007 & 0.9988 & 5 & 0.12 & 0.000 & 0.0046 & 0.0677 & 5300.47  \\
\midrule
6$\times$30 & 0 & Original & 0.0 & 0.0004 & 0.9997 & 5 & 650.04 & 0.000 & 0.0000 & 0.0000 & 650.04  \\
6$\times$30 & 0 & Original* & 6.8 & 0.0004 & 0.9997 & 5 & 136.42 & 0.000 & 0.0000 & 0.0000 & 143.23  \\
6$\times$30 & 0 & Region-reduced & 14.2 & 0.0011 & 0.9952 & 5 & 49.74 & 0.000 & 0.0014 & 0.0709 & 63.93  \\
6$\times$30 & 0 & Region-reduced*& 14.2 & 0.0011 & 0.9952 & 5 & 29.54 & 0.000 & 0.0014 & 0.0709 & 43.74  \\
6$\times$30 & 0 & Pruned & 249.0 & 0.0019 & 0.9943 & 5 & 71.47 & 0.000 & 0.0049 & 0.0446 & 320.44  \\
6$\times$30 & 0 & Pruned* & 253.4 & 0.0019 & 0.9943 & 5 & 30.54 & 0.000 & 0.0049 & 0.0446 & 283.91  \\
6$\times$30 & 0 & Student & 3131.7 & 0.0004 & 0.9997 & 5 & 62.79 & 0.000 & 0.0041 & 0.0768 & 3194.45  \\
6$\times$30 & 0 & Student* & 3136.4 & 0.0004 & 0.9997 & 5 & 29.61 & 0.000 & 0.0041 & 0.0768 & 3165.96  \\
\midrule
6$\times$30 & 0.0001 & Original & 0.0 & 0.0007 & 0.9989 & 5 & 1.44 & 0.000 & 0.0000 & 0.0000 & 1.44  \\
6$\times$30 & 0.0001 & Original* & 2.2 & 0.0007 & 0.9989 & 5 & 0.87 & 0.000 & 0.0000 & 0.0000 & 3.04  \\
6$\times$30 & 0.0001 & Region-reduced & 4.6 & 0.0011 & 0.9953 & 5 & 0.51 & 0.000 & 0.0655 & 0.1939 & 5.09  \\
6$\times$30 & 0.0001 & Region-reduced*& 4.6 & 0.0011 & 0.9953 & 5 & 0.42 & 0.000 & 0.0032 & 0.1401 & 5.00  \\
6$\times$30 & 0.0001 & Pruned & 317.4 & 0.0020 & 0.9935 & 5 & 14.56 & 0.000 & 0.0084 & 0.1388 & 332.00  \\
6$\times$30 & 0.0001 & Pruned* & 321.4 & 0.0020 & 0.9935 & 5 & 3.32 & 0.000 & 0.0084 & 0.1388 & 324.69  \\
6$\times$30 & 0.0001 & Student & 7004.8 & 0.0007 & 0.9987 & 5 & 0.13 & 0.000 & 0.0070 & 0.1658 & 7004.90  \\
6$\times$30 & 0.0001 & Student* & 7005.7 & 0.0007 & 0.9987 & 5 & 0.13 & 0.000 & 0.0066 & 0.1607 & 7005.78  \\
\midrule
8$\times$30 & 0 & Original & 0.0 & 0.0005 & 0.9996 & 0 & 7200.06 & 0.463 & 0.0000 & 0.0000 & 7200.06  \\
8$\times$30 & 0 & Original* & 16.9 & 0.0005 & 0.9996 & 3 & 2898.78 & 0.174 & -0.0006 & 0.0127 & 2915.69  \\
8$\times$30 & 0 & Region-reduced & 31.9 & 0.0009 & 0.9933 & 5 & 727.99 & 0.000 & 0.0010 & 0.0700 & 759.86  \\
8$\times$30 & 0 & Region-reduced*& 31.9 & 0.0009 & 0.9933 & 5 & 389.96 & 0.000 & 0.0010 & 0.0700 & 421.84  \\
8$\times$30 & 0 & Pruned & 310.7 & 0.0018 & 0.9952 & 5 & 1366.20 & 0.000 & 0.0081 & 0.1850 & 1676.86  \\
8$\times$30 & 0 & Pruned* & 322.2 & 0.0018 & 0.9952 & 5 & 349.84 & 0.000 & 0.0081 & 0.1850 & 671.99  \\
8$\times$30 & 0 & Student & 3140.9 & 0.0005 & 0.9996 & 4 & 2184.25 & 0.100 & 0.0038 & 0.0568 & 5325.20  \\
8$\times$30 & 0 & Student* & 3152.2 & 0.0005 & 0.9996 & 5 & 890.45 & 0.000 & 0.0022 & 0.0471 & 4042.68  \\
\midrule
8$\times$30 & 0.0001 & Original & 0.0 & 0.0008 & 0.9986 & 5 & 3.61 & 0.000 & 0.0000 & 0.0000 & 3.61  \\
8$\times$30 & 0.0001 & Original* & 3.9 & 0.0008 & 0.9986 & 5 & 1.15 & 0.000 & 0.0008 & 0.0155 & 5.03  \\
8$\times$30 & 0.0001 & Region-reduced & 8.5 & 0.0022 & 0.9838 & 5 & 0.42 & 0.000 & 0.0034 & 0.0963 & 8.90  \\
8$\times$30 & 0.0001 & Region-reduced*& 8.5 & 0.0022 & 0.9838 & 5 & 0.42 & 0.000 & 0.0034 & 0.0963 & 8.89  \\
8$\times$30 & 0.0001 & Pruned & 399.1 & 0.0023 & 0.9929 & 5 & 42.58 & 0.000 & 0.0117 & 0.0670 & 441.64  \\
8$\times$30 & 0.0001 & Pruned* & 409.0 & 0.0023 & 0.9929 & 5 & 16.86 & 0.000 & 0.0117 & 0.0670 & 425.90  \\
8$\times$30 & 0.0001 & Student & 8596.1 & 0.0009 & 0.9982 & 5 & 0.25 & 0.000 & 0.0054 & 0.0717 & 8596.37  \\
8$\times$30 & 0.0001 & Student* & 8597.4 & 0.0009 & 0.9982 & 5 & 0.16 & 0.000 & 0.0051 & 0.0693 & 8597.54  \\
\bottomrule
\end{tabular}
\end{table}

\section{Superstructure results}\label{app:superstructure}

\begin{figure}[H]
    \centering
    \includegraphics[width=.95\linewidth]{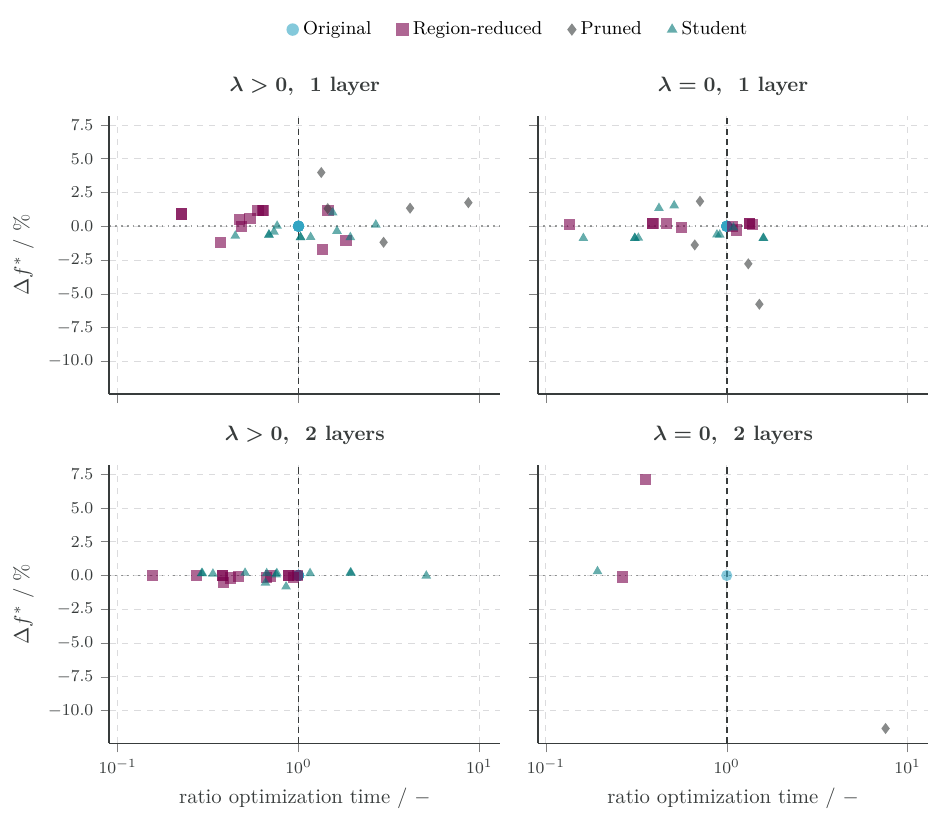}
    \caption{Relative deviation of optimal objective over ratio of optimization times obtained with the compressed models compared to the respective original models for all instances of the superstructure problem considered.}
    \label{fig:superstructure_deltaf}
\end{figure}

\begin{figure}[t]
    \centering
    \includegraphics[width=.95\linewidth]{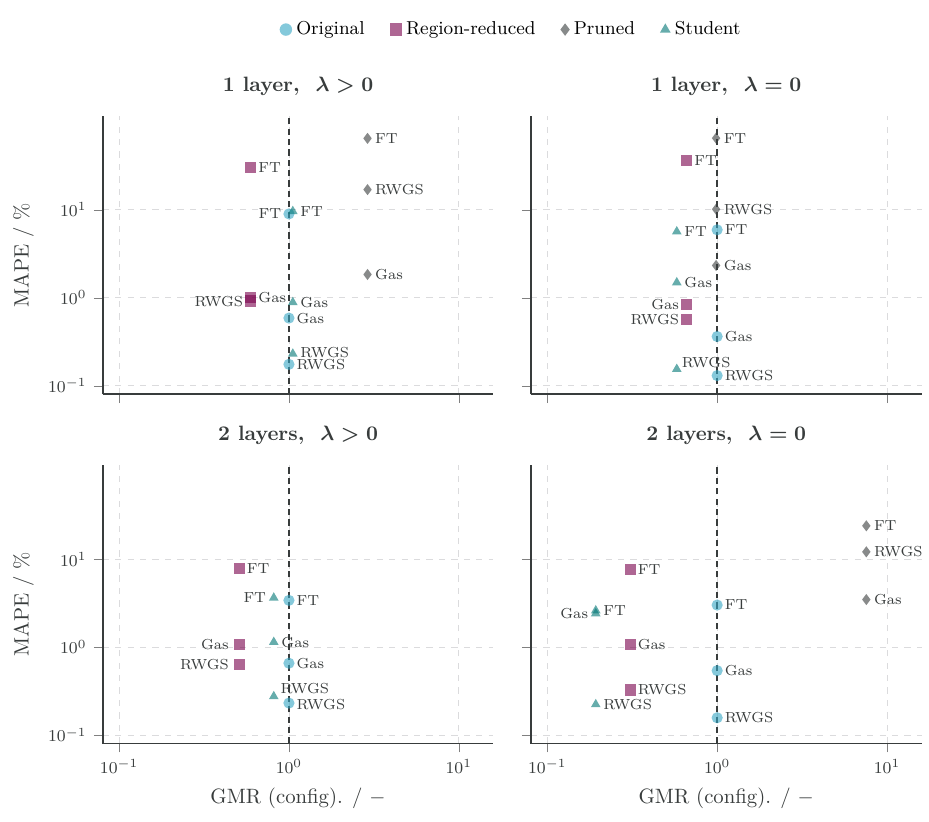}
    \caption{Mean-absolute percentage error over geometric mean of ratios of optimization times obtained with the compressed models compared to the respective original models averaged over all solved instances of superstructure problems.}
    \label{fig:superstructure_mape}
\end{figure}

For the superstructure problem, there are two ways to compare the effect of the compression methods on computational times and optimization outcomes. First, comparing against the baseline setting of \citet{Klimek2026}, i.e., ANNs with one hidden layer trained with $\ell_1$ regularization. Second, comparing against the original model for each configuration. \Cref{fig:superstructure_deltaf,fig:superstructure_mape} show the latter comparison for the superstructure results, hence the scatter plots are now centered around the original models of each configuration.

\end{appendices}

\end{document}